\documentclass[11pt]{article}
\usepackage{authblk}
\usepackage{amssymb}
\usepackage{amsmath}
\usepackage{amsfonts}
\usepackage{amscd}
\usepackage{mathtools}
\usepackage{graphicx,color}
\usepackage{fullpage}
\usepackage{tikz}
\usepackage{graphicx}
\usepackage{subcaption}
\usepackage{tabularx,ragged2e,booktabs,caption}
\usepackage{floatrow}
\usepackage{algorithm2e}

\newcommand{\norm}[1]{\left\lVert#1\right\rVert}
\newcommand{\bs}{\boldsymbol}
\newcommand{\ms}{\textsf}

\newcommand{\tb}{\textbf}
\newcommand{\ed}{{\text{d}}}

\title{A  mixed method for 3D nonlinear elasticity using finite element exterior calculus}
\author[1]{Bensingh Dhas }
\author[1]{Jamun Kumar N }
\author[1]{Debasish Roy \thanks{royd@iisc.ac.in}}
\author[2]{J N Reddy \thanks{jnreddy@tamu.edu}}
\affil[1]{CoE in Advanced mechanics of materials, Indian
Institute of Science, Bangalore 560012, India}
\affil[2]{J. Mike Walker’66 Department of Mechanical Engineering, Texas A\&M University, College Station, TX 77843-3123, USA}
\date{}

\begin{document}
\maketitle
\begin{abstract}
    This article discusses a mixed finite element technique for 3D nonlinear
    elasticity using a Hu-Washizu (HW) type variational principle. In this method, the deformed configuration and sections from its
    cotangent bundle are the input arguments. The critical points of the proposed HW
    functional enforce compatibility of these sections with the configuration,
    in addition to mechanical equilibrium and constitutive relations. This finite element (FE)
    approximation distinguishes a vector from a 1-from, a feature not commonly
    found in mixed FE approximations for nonlinear elasticity. This point of view
    permits us to construct finite elements with vastly superior performance.  Discrete approximations for the differential forms appearing
    in the variational principle are constructed with ideas borrowed from finite
    element exterior calculus. These are in turn used to
    construct a discrete approximation to our HW functional.  The discrete
    equations describing mechanical equilibrium, compatibility and constitutive
    rule, are obtained by extemizing the discrete functional with respect to
    appropriate degrees of freedom. The discrete extremization problem is
    then solved using a numerical search technique; we use
    Newton's method for this purpose. This mixed FE technique is then applied to
    a few benchmark problems wherein conventional displacement based approximations
    encounter locking and checker boarding. These studies help establish that our mixed FE approximation, which requires no artificial stabilising terms, is free of these numerical
    bottlenecks.\\
\end{abstract}
\textit{Keywords: }Hu-Washizu principle; differential forms; non-linear
elasticity; finite element exterior calculus 

\section{Introduction}
Constructing well performing numerical discretization schemes for non-linear
elasticity  is a challenging research area, having attracted enormous attention over the years
because of their usefulness in engineering design and simulation. Earlier efforts
in this area were mostly to develop finite element (FE) techniques to approximate the
deformation field alone.  However, it was soon realised that such approximations were plagued with instabilities such as shear locking, volumetric locking and checker
boarding.  Techniques based on reduced integration, enhanced and assumed strain
techniques \cite{simo1990,simo1993} were introduced to circumvent these
difficulties.  Assumed and enhanced strain techniques  introduce additional
stabilization terms in the variational formulation, the motivation for which is purely numerical. 
From the literature, it may be
gathered that constructing numerical approximations for 3-dimensional non-linear
elasticity that takes only the elasticity parameters and works well for quite a
broad range of material models, loading and applied constraints is still work in progress. This
article is an attempt, perhaps a first of its kind, to address this problem by adopting tools from geometry.
Often in mechanics, the distinction between a vector and a 1-form is understated.
This blurs the geometric meaning bought about by these
'vectors'. For example, in particle dynamics, both force and velocity may be
represented as vectors, even though this representation is physically incorrect;
it prevents us from talking about the existence of a potential for a force. If force
were to be represented as a 1-form, one could talk about its integrability in the
sense of Poincar\'e, leading to the possible existence of a potential for the force.
Moreover, the natural action of a force (a 1-form) on velocity (a vector) yields power
without the need for an inner-product. Considering force as a 1-form is common
in geometric mechanics \cite{marsden2013} and electrodynamics \cite{hehl2012}.
However, such a description is rare in elasticity and, in general, continuum
mechanics of solids.

In an effort to mimic the force-velocity pairing in particle dynamics, Frankel
\cite{frankel2011} gave a differential form based description of Cauchy and
Piola stresses (See Appendix A in \cite{frankel2011}). His idea was to decompose
stress into area and force components, where the area component was
represented by a 2-form and the force component by a 1-form. He also
suggested the use of Cartan's moving frame to describe the equation of
non-linear elasticity in arbitrary curvilinear coordinates. However this
suggestion went mostly unnoticed. Kanso \textit{et al.} also gave a description of
stress tensor using bundle valued differential forms, their description being very
similar to that of Frankel.  Building on the work of Frankel and Kanso
\textit{et al.}, Dhas and Roy \cite{dhas2020} have presented a Hu-Washizu (HW) type
mixed variational principle for nonlinear elasticity, which can accommodate a
geometric description of deformation and stresses. An important feature of this
formulation is its accommodation of the frame fields as an additional
input argument making it suitable for problems in which the affine connection
evolves, an example being Kirchhoff type shells. 

For a deformable solid, the HW functional takes variables from the base space (configuration) and fibre (tangent/ cotangent space)  of the
deformed configurations tangent bundle as independent. The HW principle then  ties these variables together with conditions which
establish compatibility between the base space and the fibre  to form a tangent bundle. These conditions often reflect certain
integrability issues, whose origin can be traced to the geometric hypothesis placed on the configuration of the body.
The simplest form of the HW variational principle is the one that takes the Jacobian of deformation as input to the
functional \cite{simo1991quasi}. This well known mixed formulation enforces the compatibility between the Jacobian
computed from deformation and the independent Jacobian field by introducing pressure as a Lagrange multiplier.
Geometrically, the Jacobian of deformation is a volume-form; the compatibility of these volume-forms require them to be
equal. In the literature, numerical techniques based on these formulations are often referred to as $p-\Theta$. It is
also well understood that constructing a numerical technique by arbitrary interpolations of $p$ and $\Theta$ leads to
instability.  This example highlights the need for a numerical method that respects the geometric structure imposed on
the fields.

Constructing finite dimensional approximations to fields that respect the
geometric structure placed on them has a long history. H Whitney
\cite{whitney2012} discussed a family of finite dimensional approximations for
differential forms on a simplicial manifold. In the modern parlance, such
approximations correspond to the lowest order FE scheme for
differential forms.  These FE spaces were developed much before the advent of
the finite element method itself. Whitney used these approximations to study the
Hodge-Laplacian.  This family of finite elements is called the
Whitey forms or Whitney elements. The geometric underpinning of this family of
finite elements went unnoticed in the literature for quite some time. Recently, 
these finite elements have come into prominence due to the effort of the
computational electromagnetics community \cite{bossavit1988}. In the finite
element literature, the FE spaces constructed by Whitey have been reinvented as the
vector finite elements of Raviart-Thomas \cite{raviart1977} and Nedelec
\cite{nedelec1980,nedelec1986}. The vector finite elements have been developed to
approximate mixed variational principles. Arnold and co-workers
\cite{arnold2010} have however unified these finite elements under the common umbrella of
finite element exterior calculus, using the language of differential forms. The last authors
interpret these finite elements as polynomial subspaces of differential forms
over $\mathbb{R}^n$. An introduction to finite element exterior calculus may be
found in \cite{ArnoldFEECBook}. Arnold and Logg \cite{feTable} have presented a
classification of finite elements using differential forms; they called this
classification the finite element periodic table. It takes into account the
geometry of the cell, the degree of the polynomial space and the degree of the
differential form. Finite elements featured in the periodic table  are complete
up to a given polynomial degree. Well known finite elements like Lagrange,
Raviart-Thomas, Nedelec, Brezzi-Douglas-Marini (BDM) etc. are all featured in the finite element periodic
table alongside the newer ones.

Apart from finite elements for differential forms, numerical approximations
based on finite difference techniques that preserve the conservation laws
exactly at the discrete level have been developed by Hyman and Shashkov
\cite{hyman1997}; see \cite{brezzi2014} for an analysis of these techniques.  In
the literature, these methods are now called mimetic techniques. Similar efforts
have been made by Hirani \cite{hiraniThesis}, who discusses a discrete version of
exterior calculus on a simplicial mesh.  In a discrete setting, Smale
\cite{smale1972} has studied  electrical circuits using ideas from algebraic
topology. Yavari \cite{yavari2008geometric} has presented a discretization for
elasticity using discrete exterior calculus. A numerical implementation of these
ideas for the equations of incompressible linearised elasticity may be found in
Angoshtari and Yavari \cite{angoshtari2013geometric}. Schrod\"er \textit{et al.}
\cite{schroder2011new} have presented a mixed finite element technique by
independently approximating deformation gradient, its cofactor and its
determinant along with the displacement field. The motivation for considering the
above variables as independent derives from the observation that these
quantities relate infinitesimal line, area and volume elements between the
reference and deformed configurations. Being mutually
dependent, their dependencies are resolved by introducing additional kinetic
variables, which enforce the compatibility.  The last authors have combined the
kinematic and kinetic variables using a HW type functional, which has all the
above kinematic and kinetic quantities as arguments. In line with the
work of Schrod\"er \textit{et al.}\cite{schroder2011new}, Bonet \textit{et al.}
\cite{bonet2016} have discussed a formulation for elastic
solids under finite deformation using a cross product defined between tensors. They claim that such an
approach would lead to simpler expressions for the first and second
derivatives of the stored energy function. Based on this algebraic construction, Bonet
\textit{et al.} \cite{bonet2015} have discussed a mixed finite element technique for
non-linear elastic solids. This mixed finite element technique utilizes the
area and volume forms obtained using the cross product definition of tensors as
addition kinematic unknowns. Conjugate generalized forces associated with these
new kinematic unknowns have been also identified.  From the perspective of the
present work, the cross product algebra of tensors introduced by Bonet
\textit{et al.} is just a specialization of the exterior algebra used in
the study of differentiable manifolds.

In a recent study, Shojaei and Yavari \cite{shojaei2019} present a mixed finite
element technique for non-linear elasticity, based on the HW variational
principle. The authors refer to their technique as the 'compatible strain mixed finite element
method (CSMFEM)'. In CSMFEM, one independently approximates deformation,
deformation gradient and Piola stress.  They base their finite element
approximation by constructing suitable finite dimensional projection operators
for a differential complex describing the equations of non-linear elasticity.
The differential complex they use is very similar to the de Rham complex
in differential topology. A difficulty with the CSMFEM is the
need for a stabilization term which penalises the incompatibility of
deformation.  Building on the idea of non-linear elasticity complex, Angoshtari
and Matin \cite{angoshtari2019} present an FE technique by identifying the
approximation spaces for deformation gradient and stress as subspaces of
$H(\text{Curl})$ and $H(\text{Div})$ respectively. However, their discrete
equilibrium, compatibility and constitutive rules cannot be obtained as
critical points of a functional. This implies that the authors are not
exploiting an important property of non-linear elasticity which is the existence
of a potential. The absence of a variational structure may pose serious difficulties
to analyse the approximation technique and place restrictions on the numerical
implementation.

In the present work, we build a finite element solution scheme for non-linear
elasticity based on a geometric description of the HW variational principle
discussed in Dhas and Roy \cite{dhas2020}. This finite element technique can be
classified as mixed, as it seeks independent approximations for both kinetic and
kinematic fields. The principal novelty of our work is that it uses a
discretization that respects the geometric information encoded in both stress
and strain fields. Instead of treating the right Cauchy-Green deformation
tensor as a single geometric object, we think of it as a composite object given
by three 1-forms. Such a point of view has the advantage of describing the
invariants of the right Cauchy-Green deformation tensor in terms of the exterior product
of the 1-forms. Similarly, the  first Piola stress may also be understood as a
composite object described by the pulled-back area-form and a traction 1-form.
Such a description of stress has never been used in computations before.

The rest of the article is organised as follows. Section 2 introduces a
geometric approach to the HW variational principle. Instead of taking
deformation gradient as the input, we take three 1-forms as input arguments.
Retrospectively, three additional 1-forms which represent the force acting on an
infinitesimal area are also taken as input. These force 1-forms enforce
compatibility between 1-forms representing deformation and the deformation field. In
Section 3, we introduce finite element exterior calculus and the finite element
spaces  $\mathcal{P}_r\Lambda^n$  and $\mathcal{P}_r^{-}\Lambda^n$ on a
3-simplex.  Using these finite elements,  we find a discrete approximation to
the HW variational principle introduced in Section 2. Newton's method is then
used to find the extremum of the discrete HW functional. The residue and tangent
operators required for Newton's method is also presented in Section 4. Numerical
simulations to validate the finite element approximation is presented in
Section 5; problems presented in this section examine the efficiency of our 
numerical technique against shear locking, volume locking and checker
boarding.  Finally, Section 6 concludes the article with a brief discussion of the
future prospects of the present work.

\section{Differential forms and Hu-Washizu variational principle}
Our present goal is to describe the HW variational principle using
differential forms. The description here will be brief; for complete details
on the kinetics and kinematics in terms of differential forms, see \cite{dhas2020}.  
We denote the reference and deformed configurations of
the body by $\mathcal{B}$ and $\mathcal{S}$ respectively. Both these
configurations can be identified with subsets of three dimensional Euclidean
space. The placements or position vectors of a material point in the reference and
deformed configurations are denoted by $X$ and  $x$ respectively. Under a
curvilinear coordinate system, coordinates of the respective position
vectors are denoted as $(X^1,...,X^3)$ and $(x^1,...,x^3)$. The deformation map
relating the position vectors of a material point in the reference and deformed
configurations is denotes by $\varphi$. In a particular coordinate system, we may
denote deformation using its components $\varphi^i$. At each tangent space of
$\mathcal{B}$ and $\mathcal{S}$, we choose a set of orthonormal vectors and
denote them as $\{E_1,...,E_3\}$ and $\{e_1,...,e_3\}$ respectively. We refer to
these sets as frames for the reference and deformed configuration. The
orthogonality noted above is with respect to the metric tensor of the Euclidean
space. The algebraic duals associated with the frame fields of the reference and deformed
configurations are denoted by $\{E^1,...,E^3\}$ and $\{e^1,...,e^3\}$; these sets
are called co-frames of the reference and deformed configurations. These frame
and co-frame fields satisfy the duality relation $E^i(E_j)=\delta^i_j$ and
$e^i(e_j)=\delta^i_j$. Using the frame and co-frame fields, differentials of
the position vector in reference and deformed configurations may be written
as,
\begin{equation}
\ed X=E_i\otimes E^i;\quad \ed x=e_i\otimes e^i.
\end{equation} 
The deformation gradient denoted by $F$ is obtained by pulling the co-vector part
of $\ed x$ back to the reference configuration, $F=e^i\otimes \varphi^*(e_i)$.
$\varphi^*(.)$ denotes the pull-back map induced by deformation. Since
we extensively employ the pull-back of deformed configuration co-frame, we reserve the
notation $\theta^i:=\varphi^*(e^i)$ and call them deformation
1-forms. In terms of deformation 1-forms, the deformation gradient can be written as,
\begin{equation}
F=e_i\otimes\theta^i.
\end{equation}
The area-forms of the reference and deformed  configurations are denoted by $
A^i$ and $a^i$ respectively. Since the dimension of these two configurations is 3,
there are three linearly independent area-forms on each. These
area-forms are given by,
\begin{align}
 A^1&=E^2\wedge E^3;\; A^2=E^3\wedge E^1;\; A^3=E^1\wedge E^2 \\
 a^1&=e^2\wedge e^3;\; a^2=e^3\wedge e^1;\; a^3=e^1\wedge e^2.
\end{align}
The symbol $(.\wedge.)$ denotes the wedge product between two differential
forms. In addition to the area-forms of the reference configuration, one may also
pull-back the area-form of the deformed configuration to the reference
configuration, $ \ms{A}^i=\varphi^*(a^i)$. These pulled-back area-forms
may be written as,
\begin{equation}
\ms{A}^1=\theta^2\wedge\theta^3;\;
\ms{A}^2=\theta^3\wedge\theta^1;\;
\ms{A}^3=\theta^1\wedge\theta^2.
\end{equation}
Note that, the pulled-back area-forms span the same linear space spanned by
$A^i$. Volume-forms on the reference and deformed configurations are denoted
by $V$ and $v$ respectively. Similar to the area-forms, one may also pull the volume-form of the deformed
configuration back to the reference; we denote the
pulled-back volume-form by $ \ms{V}=\varphi^*(v)$. In terms of deformation
1-forms, $ \ms{V}$ may be written as,
\begin{equation}
 \ms{V}=\theta^1\wedge\theta^2\wedge\theta^3.
\end{equation}
In terms of deformation 1-forms and pulled-back area forms, the invariants of
the right Cauchy-Green deformation tensor may be written as,
\begin{align}
    \label{eq:I1}
I_1&=\left\langle\theta^i,\theta^i \right\rangle\\
    \label{eq:I2}
I_2&=\left\langle  \ms{A}^i, \ms{A}^i \right\rangle\\
    \label{eq:I3}
I_3&={}_{\star}(\theta^1\wedge\theta^2\wedge\theta^3)
\end{align}
In \eqref{eq:I3}, ${}_{\star}(.)$ denotes the Hodge star operator, which defines
an isometry between the spaces of $n$ and $m-n$ differential forms.  In defining $I_1$ and
$I_2$, we have used the inner-product induced by the Euclidean metric on the
spaces of 1- and 2-forms.

\subsection{Hu-Washizu variational principle}
In its original form, the HW functional for a hyperelastic elastic solid takes the
deformation field, deformation gradient and first Piola stress as inputs. 
However, in keeping with our setting, we adopt a HW type functional whose input arguments are
deformation field and differential forms describing deformation and stress. In
terms of these variables, the HW functional may be written as,  
\begin{equation}    
    \Pi(\theta^i,t^i,\varphi)=\int_{\mathcal{B}}
    W(\theta^{i}) \ed V  - ( t^{i}\otimes \mathsf{A}^{i}) \dot{
        \wedge}(e_{j}\otimes(\theta^{j}-\ed\varphi^{j}))+\int_{\partial
    \mathcal{B}_t}\left\langle t^\sharp,\varphi \right\rangle
    \label{eq:HWprinciple}
\end{equation}
In the above equation, $t^i$ denotes the traction 1-form. Instead of working with the
conventional definition of stress, we use a geometric approach to stress originally
put fourth by Frankel (see Appendix-C of the book) \cite{frankel2011} and later
explored by  Kanso \textit{et al.} \cite{kanso2007}. The idea behind this
construction is to decompose the stress tensor into traction and area forms.
This construction is also intuitive, in line with the understanding of stress as
traction (force) acting on an infinitesimal area. Using this, the
first Piola stress may be defined as the partial pull-back of the Cauchy stress,
where the pull-back is applied only to the area leg, leaving the traction 
undisturbed (for more details see \cite{mfe} and \cite{kanso2007}).  Thus the
first Piola stress is related to the traction 1-from and pulled-back area-from
in the following way,
\begin{align}
    \nonumber
    P&=t^i\otimes \varphi^*(a^i)\\
   \label{eq:piolaTransform}
     &=t^i\otimes \ms{A}^i.
\end{align}
In \eqref{eq:HWprinciple}, $\partial \mathcal{B}_t$ denotes the traction
boundary on with $t$ is applied and $(.)^\sharp$ denotes the sharp operator
introduced to make the contraction consistent. In writing the functional given
in \eqref{eq:HWprinciple}, we have assumed that the response of the body is described
by a stored energy density function $W$. Conventionally, $W$ 
depends on the right Cauchy-Green deformation tensor. In the present work, we
assume the stored energy function to be a function of $\theta^i$. However, the
dependence is only through the invariants of the right Cauchy-Green deformation
tensor described in \eqref{eq:I1}, \eqref{eq:I2} and \eqref{eq:I3}. The symbol
$\dot{\wedge}$ used in \eqref{eq:HWprinciple} denotes the bilinear product,
whose definition is somewhat similar to the one given in by Kanso \textit{et
al.}\cite{kanso2007}.  This algebraic operation is useful in writing the work
done by the Piola stress on deformation. For differential forms $a,b,c$ and
the vector field $v$, we define $\dot{\wedge}$ as $ (a\otimes b)
\dot{\wedge}(v\otimes c)=a(v)(b \wedge c)$. Our definition of $\dot{\wedge}$ does
not rely on the metric tensor. The definition of $\dot{\wedge}$ given by  Kanso
\textit{et al.} uses the metric tensor, which implies that work done by 
stress is dependent on the metric structure placed on the configuration
manifold. However, from classical dynamics, it is well known that one can compute
power or work without the metric tensor.

Instead of postulating the relationship between these differential forms, we
follow a variational route to obtain these relations, which 
are obtained as the critical points of the HW functional. The critical points of
the HW energy functional \eqref{eq:HWprinciple} are thus given as, 
\begin{align}
    \label{eq:compatibility}
    \delta_{t^i}\Pi&: \theta^i-\ed \varphi^i=0;\quad i=1,...,3\\
    \label{eq:constitutive1}
    \delta_{\theta^1}\Pi&: \frac{\partial W}{\partial \theta^1}
        -[t^1(e_1){}_{\star}(\theta^2\wedge\theta^3)
        +t^2(e_1){}_{\star}( \theta^3\wedge\theta^1)
        +t^3(e_1){}_{\star}(\theta^2\wedge\theta^1)
        ]^\sharp=0\\
    \label{eq:constitutive2}
     \delta_{\theta^2}\Pi&:   \frac{\partial W}{\partial \theta^2}
        -[t^1(e_2){}_{\star}(\theta^2\wedge \theta^3)
        +t^2(e_2){}_{\star}((\theta^3\wedge \theta^1)
        +t^3(e_2){}_{\star}(\theta^1\wedge \theta^2)]^\sharp=0\\
    \label{eq:constitutive3}
      \delta_{\theta^3}\Pi&:  \frac{\partial W}{\partial \theta^3}
        -[t^1(e_3){}_{\star}(\theta^2\wedge \theta^3) 
        +t^2(e_3){}_{\star}(\theta^3\wedge \theta^1)
        +t^3(e_3){}_{\star}(\theta^1\wedge\theta^2) ]^\sharp=0\\
        \label{eq:equilibrium}
    \delta_{\varphi^i}\Pi&:\ed (t^j(e_k))\wedge \mathsf{A}^j=0;\quad
    i=1,...,3.
\end{align}    
In the above equation, $\delta_{(.)}W$ denotes the Gateaux derivative of $W$ in the
direction of an input argument. Eq. \eqref{eq:compatibility} describes the compatibility
between the deformation 1-forms and deformation, Eqs. \eqref{eq:constitutive1}
\eqref{eq:constitutive2} and \eqref{eq:constitutive3} are the constitutive rules
relating the deformation 1-forms and traction 1-forms and Eq. \eqref{eq:equilibrium}
describes mechanical equilibrium.  For more details on the variational
principle and calculation of critical points of the HW functional, the reader may see \cite{dhas2020}.

\section{Finite element approximation of differential forms}
In this section, we discuss a finite dimensional approximation for differential
forms using finite element spaces; these approximation spaces can be neatly
dealt within the classical definition of finite elements provided by Cairlet
\cite{ciarletbook}. Recall that a finite element is a triplet denoted by
$(K,\mathcal{P}(K),\Sigma )$, where $K$  denotes a triangulation or a cell
approximation of a given domain, $\mathcal{P}(K)$ denotes the space of piecewise
polynomials of degree less than or equal to $n$ defined on $K$ and $\Sigma$
denotes the degrees of freedom, whose values are from the dual space of
$\mathcal{P}(K)$.  By choosing the degrees of freedom, one prescribes how the
piecewise polynomial spaces are glued together globally. The choice of degrees
of freedom also provides a basis for the finite element space.  Lagrange,
Raviart-Thomas and N\'edelec finite elements are common examples of finite
elements that may be discussed using the definition given
above. In this work, we use a simplicial approximation for the domain. The
finite element spaces are constructed on a reference simplex and then mapped
back to the actual simplex using an affine map. Notations used in this
section closely follow those developed by Arnold and co-workers
\cite{arnold2010,ArnoldFEECBook}. 

\subsection{Spaces $\mathcal{P}_r\Lambda^n$ and $\mathcal{P}_r^{-}\Lambda^n$}
We denote the space of $m$-variable polynomials of degree $r$  in $\mathbb{R}^m$
by $\mathcal{P}_r(\mathbb{R}^m)$. The space of polynomial differential forms
with form degree $n$ and polynomial degree $r$ is denoted by
$\mathcal{P}_r\Lambda^n(\mathbb{R}^m)$. Often we suppress $\mathbb{R}^m$ from
our notation and simply denote these spaces by $\mathcal{P}_r$ and
$\mathcal{P}_r\Lambda^n$ respectively.  By a polynomial differential form, we mean
the following: for $v_1,...,v_n \in$ $T\mathbb{R}^m$,
\begin{equation}
    \mathcal{P}_r\Lambda^n=\{\omega \in \Lambda^n(\mathbb{R}^m)| \omega(v_1,...,v_n)\in
    \mathcal{P}_r\}.
\end{equation}
In other words, the coefficient functions of $\omega$ are polynomials of degree
$r$. The space $\mathcal{P}_r$ has dimension 
$\binom{r+m}{m}$. The space
$\mathcal{P}_r\Lambda^n$ can be constructed as a product of $m-$ variable
polynomials and $n-$degree differential forms. From this construction, the
dimension of $\mathcal{P}_r\Lambda^n$ can be computed as $\binom{m+r}{m}\binom{n}{m}$.

We denote the Koszul operator by $\kappa$; the action of $\kappa$ on a
differential form $\omega$ is given by,
\begin{equation}
    \kappa \omega=\omega(v_X,v_1,...,v_{m-1})
\end{equation}
where, $v_X$ denotes the vector which translates the point $X\in\mathbb{R}^m$ to
the origin.  Note that $\kappa$ decreases the degree of a differential form by
one; in other words, the operator $\kappa$ has degree $-1$. It also increases
the polynomial degree of $\omega$ by one. It is readily verifiable that the operator
$\kappa$ commutes with affine pull-back; if $L:\mathbb{R}^n\rightarrow
\mathbb{R}^n$ is an affine linear map,  then $L^{*}\kappa \omega = \kappa L^{*}
\omega$. 

Another important subspace of polynomial differential forms discussed by Arnold
and co-workers is $\mathcal{P}_r^-\Lambda^n$. The dimension of this space is
larger than $\mathcal{P}_{r-1}\Lambda^n$, but smaller than
$\mathcal{P}_{r}\Lambda^n$. The  space $\mathcal{P}_r^{-}\Lambda^n$ is defined
as,
\begin{equation}
    \mathcal{P}_r^{-}\Lambda^n=\{\omega \in \mathcal{P}_r\Lambda^k|\kappa \omega\in
    \mathcal{P}_{r}\Lambda^{n-1}\}.
\end{equation}
The dimension of the space $\mathcal{P}_r^-\Lambda(\mathbb{R}^m)$ can be
computed as $\binom{r+k-1}{k}\binom{m+r}{m-k}$.
For differential form of
degrees 0 and $m$, we have $\mathcal{P}_r^{-}\Lambda^0=\mathcal{P}_r\Lambda^0$
and $\mathcal{P}_r^{-}\Lambda^m=\mathcal{P}_r\Lambda^m$. The spaces
$\mathcal{P}_r\Lambda^n$ and $\mathcal{P}_r^-\Lambda^n$ are both affine invariant.
When restricted to a simplex of dimension $m$, these spaces contain well known
finite element spaces as special cases for the specific choice of polynomial degree
and form degree.

\subsection{Simplicial approximations for $\mathcal{B}$ and $\mathcal{S}$}
An $n$-simplex in $\mathbb{R}^m$, $m\geq n$ is the closed convex-hull generated
by $n+1$ points (vertices) in $\mathbb{R}^m$. We may also refer to an $n-$simplex by a
set containing $n+1$ integers indexing the $n+1$ points.  A simplicial complex
$K$ of dimension $n$ is a collection of $n$ simplicies such that every face of
$K$ is also a simplex in $K$ and if two simplices in $K$ intersect, the
intersection is also a simplex in $K$. A three dimensional finite element mesh,
constructed using tetrahedra is an example of simplicial complex. The
simplicial approximation for the reference and deformed configurations are
denoted by $K_\mathcal{B}$ and $K_\mathcal{S}$ respectively. The finite element
basis functions are constructed on a reference tetrahedron, which we denote by $T_R$. 
These basis functions are then pulled-back to
the actual tetrahedron using an affine map. The barycentric coordinates of $T_R$
are denoted by $\lambda^i$, $0\leq i \leq 3$. These coordinates satisfy the well known 
constraint $\sum_{i=0}^n\lambda^i=1$. Fig. \ref{fig:orientationTet}
shows the reference tetrahedron and the numbering used to index the $0-$ and
$1-$ simplices.  Table \ref{table:dataSimplices} lists the coordinates of
the 0-simplex and the vertices that make up the 1-simplices.  

\begin{center}
    \captionof{table}{Sub-simplices of dimension zero and one for a $3-$
    simplex(tetrahedron)}
 \begin{tabular}{cccccc} \hline
     Vertex index & coordinate & Edge index & [$v_1$,$v_2$] \\ \hline
      1 & (0,0,0) &1 & [1,2] \\
      2 & (1,0,0)  &2 &[1,3] \\
      3 & (0,1,0)  &3 &[1,4] \\
      4 & (0,0,1)  &4 &[2,3] \\
        & &5 &[2,4] \\
        & &6 &[3,4] \\
      \hline
  \end{tabular}
    \label{table:dataSimplices}
\end{center}

\subsection{Orientation of simplices}
\begin{figure}
    \centering
    \includegraphics[scale=0.6]{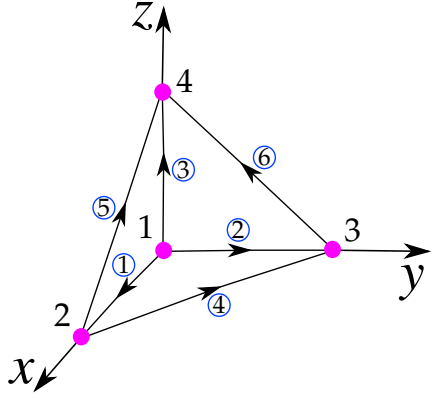} 
    \caption{Reference tetrahedron used to construct the basis for the finite
    element spaces. Edge numbers are marked with circles, while an arrow along  an
    edge indicates a positive orientation.}
    \label{fig:orientationTet}
         
\end{figure}

The idea of orientation is crucial in defining the global degrees of freedom for a
finite element differential form. We say that a manifold  is orientable if it admits
a global volume form. In the present context, we define a simplex to have a
positive orientation if the order of the vertices follow the lexicographic order
or an even permutation of the lexicographic order.  Using this definition, it follows 
that there are only two orientations possible for a simplex of any
dimension. However, the scenario is quite different when one turns to a simplicial
complex constructed using simplices. Orientation is a global property; it may
turn out that a simplicial complex may not admit an orientation. We refer to
such complexes as non-orientable.  A triangulation placed on a M\"obious strip
is an example of a simplicial complex which does not admit an orientation. We
exclude these non-orientable simplicial complexes from further discussions.
Figure \ref{fig:orientationTet} shows the numbering scheme we use for a
3-simplex; an arrow placed along an edge indicates positive orientation.

\subsection{Bases for space $\mathcal{P}_r\Lambda^n(T_R) $ and
$\mathcal{P}_r^{-}\Lambda^n(T_R)$}
In the previous subsections, we have introduced the simplicial approximation for a
configuration and the spaces $\mathcal{P}_r\Lambda^n$ and
$\mathcal{P}_r^{-}\Lambda^n$. Presently, we discuss finite element
basis functions for the spaces  $\mathcal{P}_r\Lambda^n(T_R) $ and
$\mathcal{P}_r^{-}\Lambda^n(T_R)$. For differential forms $\omega, \eta\in
\Lambda^n$, the $L^2$ inner product between them is defined as,
\begin{eqnarray}
    \left\langle\omega,\eta\right\rangle_{L^2}=\int\omega\wedge{}_{*}\eta
\end{eqnarray}
The space $L^2\Lambda^n$ is obtained by completing the space $\Lambda^n$ in the
norm induced by the inner product $\left\langle .,.\right\rangle_{L^2}$. We
denote the Hilbert space of differential forms with form degree $n$ by
$H\Lambda^n$, which is defined as,
\begin{equation}
    H\Lambda^n=\{\omega\in L^2 \Lambda^{n}|\ed \omega \in L^2 \Lambda^{n+1}\}.
\end{equation}
The polynomial spaces $\mathcal{P}_r\Lambda^n(K)$ and
$\mathcal{P}_r^{-}\Lambda^n(K)$ are constructed by restricting the spaces
$\mathcal{P}_r\Lambda^n(\mathbb R^m)$ and $\mathcal{P}_r^{-}\Lambda^n(\mathbb
R^m)$ to $K$. However, within a finite element framework, the spaces
$\mathcal{P}_r^-{\Lambda}^n(K)$ and $\mathcal{P}_r{\Lambda}^n(K)$ are obtained
by gluing together polynomial spaces defined on each simplex $T\in K$. Since we
use an affine map to pull the finite element basis from the reference simplex,
we need to make sure that the pulled back differential forms are elements in
$H\Lambda^n$. The following condition is used to enforce this requirement. If
$\omega \in H\Lambda^n(K)$, $T_1, T_2$ are simplices in $K$ with a common face
$f$ of dimension $m-1$, then $\omega|_{f\in T_1}=\omega|_{f\in T_2}$.

Arnold \textit{ et al.} \cite{arnold2009geometric} proposed a geometric
decomposition for the dual space of the polynomial spaces
$\mathcal{P}_r\Lambda^n$ and $\mathcal{P}_r^-\Lambda^n$ by mimicking the
geometric decomposition for the degrees of freedom for the Lagrange finite
element space, leading to the Bernstein basis. For a polynomial differential
form with a given degree, they distributed the degrees of freedom on all
sub-simplices with dimension greater than or equal to the form degree. This
decomposition of the degrees of freedom affords a convenient computational
basis index by the sub-simplices. Table \ref{table:feSpaceAndDim} gives the
basis functions for finite element differential forms of degree 0 and 1, in
terms of the barycentric coordinates of the 3-simplex.

\begin{minipage}{\linewidth}
\centering
\captionof{table}{FE basis functions with form degree 0 and 1 and polynomial
    degree 1} \label{tab:basis1Form} 
\begin{tabular}{c c l l }
\cline{1-4}
    FE spaces & \# DoF & Node [i]&Edge $[i, j]$\\ 
\cline{1-4}
    $\mathcal{P}_{1}\Lambda^{0}(\mathcal{T})$ & 4 & $\lambda^i$ &-\\
    $\mathcal{P}^{-}_{1}\Lambda^{1}(\mathcal{T})$      & 6   & -&$
        \lambda^{i} \ed\lambda^{j} - \lambda^{j}\ed \lambda^{i}$\\ 
    $\mathcal{P}_{1}\Lambda^{1}(\mathcal{T})$  & 12  &-& $
        \lambda^{i}\ed\lambda^{j} , \lambda^{j}\ed \lambda^{i}$\\ 
\bottomrule[.25pt]
\end{tabular}
    \label{table:feSpaceAndDim}
\end{minipage}
It should be mentioned that the above basis functions do not take into account
the orientation of the simplex; multiplying the basis function by a -1 corrects
it for orientation.

\section{Discretization of HW variational principle}
We now construct a discrete approximation to the HW functional
developed in Section 2. To do this we replace the input arguments of  the HW
functional by their discrete approximation. Given a simplicial approximation for
the reference configuration of an elastic body, our objective is to find a
simplicial approximation for the deformed configuration, such that the
approximation extremizes the discrete HW functional.  In addition to this, one
also needs to find  approximations for the traction and deformation 1-forms. In
this work, we use finite element spaces with polynomial degree 1 to approximate
deformation and traction 1-forms; displacements are also approximated by
polynomials of the same degree. In Section 2, we interpreted traction 1-form
as a Lagrange multiplier enforcing the constraint between deformation and
deformation 1-forms. This understanding guides us to choose the finite
dimensional approximation for traction 1-forms to come from an FE space that is
smaller than the FE space for deformation 1-forms. We use the finite element space
$\mathcal{P}_1\Lambda^1$ to interpolate deformation 1-forms and
$\mathcal{P}_1^{-}\Lambda^1$ to approximate traction 1-forms; the dimension of
these spaces may be found in \ref{table:feSpaceAndDim}. In our present theory,
each displacement component is understood as a zero-form and hence approximated
using the FE space $\mathcal{P}_1\Lambda^0$. Table \ref{table:fieldsDof} gives a summary
of the FE spaces used to approximate different fields involved in the HW
functional.
\begin{minipage}{\linewidth}
\centering
\captionof{table}{Fields and FE spaces used to
    approximate them} 
\begin{tabular}{l c c }
\cline{1-3}
Field & FE space & \# Dof \\ 
\cline{1-3}
    Displacement& $\mathcal{P}_{1}\Lambda^{0}$ & 4 \\
    Traction 1-form& $\mathcal{P}^{-}_{1}\Lambda^{1}$ & 6  \\ 
    Deformation 1-form& $\mathcal{P}_{1}\Lambda^{1}$  & 12 \\ 
\bottomrule[.25pt]
\end{tabular}
    \label{table:fieldsDof}
\end{minipage}
We introduce the following notations so that the calculations to follow are more
readable. The basis functions from the spaces $\mathcal{P}_1^{-}\Lambda^1$  and
$\mathcal{P}_1\Lambda^1$ are denoted by $\phi^i$ and $\psi^i$ respectively. Using
this new notation, the discrete approximation for different fields involved in
the HW functional may be written as,
\begin{align}
    \label{eq:approxTheta}
    \theta^i_h&:=\sum_{j=1}^{12}\phi^j \bar{\theta}_j^i;\quad \phi^j\in
    \mathcal{P}_1\Lambda^1(T_{\mathcal{B}}),\\
    \label{eq:approxTraction}
    t^i_h&:=\sum_{j=1}^{6}\psi^j \bar{t}_j^i;\quad \psi^j\in
    \mathcal{P}_1^{-}\Lambda^1(T_{\mathcal{B}}),\\
    \label{eq:approxU}
    u^i_h&:=\sum_{j=1}^{4}\lambda^j \bar{u}^i_j ;\quad \lambda^j\in
    \mathcal{P}_1\Lambda^0(T_{\mathcal{B}}).
\end{align}
In the above equations, the subscript $h$ is introduced to indicate that the right
hand side only approximates the actual field. 
$\bar{\theta}^i$, $\bar{t}$ and $\bar{u}^i$ denote the degrees of freedom
associated with the $i$th deformation and traction 1-form and displacement field
respectively. The differential to the displacement approximation can be written
as,
\begin{equation}
    \ed u^i_h=\sum_{j=1}^{4} \ed\lambda^j\bar{u}^i_j.
    \label{eq:approxdu}
\end{equation}
From Table \ref{table:fieldsDof}, one may evaluate our
present approximation to have 66 DoFs on a 3-simplex. Out of these, 54 are
identified with the edges (1 dimensional subsimplices) and the remaining 12 with
vertices (0 dimensional subsimplices). In the matrix notation,
\eqref{eq:approxTheta}, \eqref{eq:approxTraction}, \eqref{eq:approxU} and  \eqref{eq:approxdu} may be
written as,
\begin{equation}
    \theta^i_h=\bs{\phi}\bs{\theta}^i;\quad t^i_h=\bs{\psi}\tb{t}^i;\quad
    u^i_h=\tb{N}\tb{u}^i;\quad \ed u^i_h=\ed \tb{N}\tb{u}^i.
\end{equation}
where, $\bs{\phi}$, $\bs{\psi}$, $\tb{N}$ and $\ed\tb{N}$ are matrices with
sizes 3 x 12, 3 x 6, 1 x 4 and 3 x 4 respectively. These matrices have the shape
functions of the associated FE space as their columns, i.e. $\bs{\phi}$ has
$\phi^i$ as its $i$th column and so on. These matrices act on the DoF vector to produce the discrete approximation for the field. The DoF
vectors $\bs{\theta}^i$, $\tb{t}^i$ and $\tb{u}^i$ are obtained by stacking the
DoFs of the $i$th field one over the other.  Using these approximations, the
discrete  HW functional may be written as,
\begin{equation}
    I^h=\int_{\mathcal{T}_{\mathcal{B}}}W^h \ed V-(t^i_h\otimes\ms
    A^i_h)\dot{\wedge}(e_j \otimes (\theta^j_h-\ed \varphi^j_h))+\int_{(\partial
    \mathcal{T}_{\mathcal{B}})_{t}} \left\langle t^{\sharp},\varphi^h\right\rangle.
    \label{eq:discreteHW}
\end{equation}
As usual, we have introduced  $h$ as a superscript to indicate that the above
functional is only an approximation to the original. Throughout the
computations, we assume that the frames $e_i$ are fixed and rectilinear.
\footnote{A frame is rectilinear if its orientation does not change across the
manifold.} 

\subsection{Newton's method}
Given a simplicial approximation for the reference configuration of the elastic
body, we use Newton's method to extremize the discrete HW variational
functional. The discrete approximation to displacement determines a
simplicial approximation for the deformed configuration, while those for 
deformation and traction 1-forms determine
stress and strain states of the body. Newton's method is a quadratic
search technique requiring the first and the second derivatives of the
objective functional. In computational mechanics literature, these derivatives are
often called residue vector and tangent operator. We now compute the first
and second derivatives of the discrete HW functional using the notion of Gateaux
derivative. Since we are working in a finite dimensional setting, the Gateaux
derivative is much simpler than the one discussed in the previous section.
To highlight this difference, we denote the finite dimensional version of the
Gateaux derivative by $\ms D_{(.)}I^h$. Similarly the second derivatives are
denoted by $\ms D_{(.)(.)}I^h$, where $(.)$ denotes the DoF vector with respect
to which the derivative is calculated. The scheme we use to extremize the
discrete HW functional is presented in Algorithm \ref{algo:NewtonsMethod}. At
each Newton iteration,  incremental DoFs computed are denoted by prefixing a $\Delta$
before the respecting DoF. These incremental DoFs are obtained by solving the system of linear
equations $\mathcal{K}X=-\mathcal{R}$, where $X$ denotes the incremental
DoF vector given by, $X=[\Delta \bs{\theta}^1,\Delta \bs{\theta}^2,\Delta
\bs{\theta}^3,\Delta \tb{t}^1,\Delta \tb{t}^2,\Delta \tb{t}^3,\Delta
\tb{u}^1,\Delta \tb{u}^2,\Delta \tb{u}^3]^t$. The components of the residue
vector may be formally written as,
\begin{equation}
    \mathcal{R}:=[\ms{D}_{\bs \theta^1} I^h, \;
    \ms{D}_{\bs \theta^2} I^h, \;
    \ms{D}_{\bs \theta^3} I^h, \;   
    \ms{D}_{\tb t^1}      I^h, \; 
    \ms{D}_{\tb t^2}      I^h, \; 
    \ms{D}_{\tb t^3}      I^h, \; 
    \ms{D}_{\tb u^1} I^h,\;
    \ms{D}_{\tb u^2} I^h,\;
    \ms{D}_{\tb u^3} I^h]^t.
    \label{eq:residue}
\end{equation}

\begin{algorithm}[H]
	\SetAlgoLined
	 \KwData{\\
      1. Simplicial approximation for $\mathcal{B}$\;
	  2. Displacement $\tb{u}^i_{prev}$ from the previous load-step\;
	  3. Traction 1-form $\tb{t}^i_{prev}$ from the previous load-step\;
      4. Deformation 1-form $\bs{\theta}^i_{prev}$ from the previous load-step\;
      5. Displacement and traction boundary conditions for the current load step\;
	}
	 \KwResult{$\bs{\theta}^i, \tb{t}^i$ and $\tb{u}^i$ for the prescribed displacement and traction conditions}
	 Initialize  $\mathcal{K}$ as a sparse matrix and $\mathcal{R}$ as a sparse vector\;
     Initialize  $TOL$\;
     $\tb{u}^i\leftarrow\tb{u}^i_{prev}$\;
     $\tb{t}^i\leftarrow\tb{t}^i_{prev}$\;
     $\bs{\theta}^i\leftarrow\bs{\theta}^i_{prev}$\;

    Modify $\mathcal{R}$ to incorporate traction and displacement conditions\;

	 \While{$\norm{\mathcal{R}} \geq TOL$}{
		\tb{Initialization:} $itEl=1$\\	
		
		 \While{itEl $\leq$ \# Elements}{
                Compute element $K$\ using \eqref{eq:tangent}\;
                Compute element $R$ using \eqref{eq:residue}\;
                Assemble $\mathcal{K}:\mathcal{K}\leftarrow K$\;
                Assemble $\mathcal{R}:\mathcal{R}\leftarrow R$\;
               }
               Apply displacement boundary conditions to $\mathcal{K}$\;
               Solve: $\mathcal{K}X=-\mathcal{R}$\;
               $\bs{\theta}^i\leftarrow\bs{\theta}^i+\Delta\bs{\theta}^i$\;
               $\tb{t}^i\leftarrow\tb{t}^i+\Delta\tb{t}^i$\;
               $\tb{u}^i\leftarrow\bs{u}^i+\Delta\bs{u}^i$\;
	  
	 }
         \caption{Newton's method for extremising the discrete HW functional}
         \label{algo:NewtonsMethod}
	\end{algorithm}

    
    

The first derivative of the discrete HW functional with respect to the
DoFs associated with the deformation 1-forms can be computed as,
\begin{align}
    \nonumber
    \ms{D}_{\bs{\theta}^1} I^h&=\int_{\mathcal{T}_{\mathcal{B}}}\ms{D}_{\bs{\theta}^1}W^h
    -(t^1(e_1)+t^2(e_2)+t^3(e_3))\bs{\phi}\wedge\theta^2\wedge\theta^3\\
    \nonumber
    &-t^2(e_1)(\bs \phi\wedge\theta^3\wedge\ed \varphi^1)
    -t^2(e_2)(\bs \phi \wedge\theta^3\wedge\ed \varphi^2)
    -t^2(e_3)(\bs \phi \wedge \theta^3 \wedge \ed \varphi^3)\\
    &+t^3(e_1)(\bs \phi\wedge\theta^2\wedge\ed \varphi^1)
    +t^3(e_2)(\bs \phi \wedge\theta^2\wedge\ed \varphi^2)
    +t^3(e_3)(\bs \phi \wedge \theta^2 \wedge \ed \varphi^3),\\
    \nonumber
    \ms{D}_{\bs{\theta}^2}I^h&=\int_{\mathcal{T}_{\mathcal{B}}}\ms{D}_{\bs{\theta}^2}W^h
    +(t^1(e_1)+t^2(e_2)+t^3(e_3))\bs{\phi}\wedge\theta^1\wedge\theta^3\\
    \nonumber
    &+t^1(e_1)(\bs \phi\wedge\theta^3\wedge\ed \varphi^1)
    +t^1(e_2)(\bs \phi \wedge\theta^3\wedge\ed \varphi^2)
    +t^1(e_3)(\bs \phi \wedge \theta^3 \wedge \ed \varphi^3)\\
    &-t^3(e_1)(\bs \phi\wedge\theta^1\wedge\ed \varphi^1)
    -t^3(e_2)(\bs \phi \wedge\theta^1\wedge\ed \varphi^2)
    -t^3(e_3)(\bs \phi \wedge \theta^1 \wedge \ed \varphi^3),\\
    \nonumber
    \ms{D}_{\bs{\theta}^3}I^h&=\int_{\mathcal{T}_{\mathcal{B}}}\ms{D}_{\bs{\theta}^3}W^h
    -(t^1(e_1)+t^2(e_2)+t^3(e_3))\bs{\phi}\wedge\theta^1\wedge\theta^2\\
    \nonumber
    &-t^1(e_1)(\bs \phi\wedge\theta^2\wedge\ed \varphi^1)
    -t^1(e_2)(\bs \phi \wedge\theta^2\wedge\ed \varphi^2)
    -t^1(e_3)(\bs \phi \wedge \theta^2 \wedge \ed \varphi^3)\\
    &+t^2(e_1)(\bs \phi\wedge\theta^1\wedge\ed \varphi^1)
    +t^2(e_2)(\bs \phi \wedge\theta^1\wedge\ed \varphi^2)
    +t^2(e_3)(\bs \phi \wedge \theta^1 \wedge \ed \varphi^3).
\end{align}
In the above equations, $\bs\phi\wedge \theta^i \wedge
\theta^j:=\sum_{k=1}^{12}\phi^k\wedge\theta^i\wedge\theta^j$. Similarly,
$\bs\phi\wedge \theta^i \wedge \ed
\varphi^j:=\sum_{k=1}^{12}\phi^k\wedge\theta^i\wedge\ed \varphi^j$. The first
derivatives of the discrete HW functional with respect to the DoFs assocaitated
with the traction 1-form are given as,
\begin{align}
    \ms D_{\tb{t}^1}I^h&=\int_{\mathcal{T}_{\mathcal{B}}}
    -(\bs \psi e_1)(\theta^1\wedge\theta^2\wedge\theta^3)
    +(\bs \psi e_1)(\theta^2\wedge\theta^3\wedge\ed \varphi^1)
    +(\bs \psi e_2) (\theta^2\wedge\theta^3\wedge\ed \varphi^2)
    +(\bs \psi e_3) (\theta^2\wedge\theta^3\wedge\ed \varphi^3),\\
    \ms D_{\tb{t}^2}I^h&=\int_{\mathcal{T}_{\mathcal{B}}}
    -(\bs \psi e_2)(\theta^1\wedge\theta^2\wedge\theta^3)
    +(\bs \psi e_1)(\theta^3\wedge\theta^1\wedge\ed \varphi^1)
    +(\bs \psi e_2) (\theta^3\wedge\theta^1\wedge\ed \varphi^2)
    +(\bs \psi e_3) (\theta^3\wedge\theta^1\wedge\ed \varphi^3),\\
    \ms D_{\tb{t}^3}I^h&=\int_{\mathcal{T}_{\mathcal{B}}}
    -(\bs \psi e_3)(\theta^1\wedge\theta^2\wedge\theta^3)
    +(\bs \psi e_1)(\theta^1\wedge\theta^2\wedge\ed \varphi^1)
    +(\bs \psi e_2) (\theta^1\wedge\theta^2\wedge\ed \varphi^2)
    +(\bs \psi e_3) (\theta^1\wedge\theta^2\wedge\ed \varphi^3).
\end{align}
In the above equations, $\bs \psi e_i$ is a vector whose components are
$\psi^i(e_i)$. The first derivatives of $I^h$ with respect to the displacement
DoFs may be computed as,
\begin{align}
    \ms D_{\tb u^1}I^h&=\int_{\mathcal{T}_{\mathcal{B}}}
    t^1(e_1)(\ed \tb N\wedge\theta^2\wedge\theta^3)
    +t^2(e_1)(\ed \tb N\wedge\theta^3\wedge\theta^1)
    +t^3(e_1)(\ed \tb N\wedge\theta^1\wedge\theta^2),\\
    \ms D_{\tb u^2}I^h&=\int_{\mathcal{T}_{\mathcal{B}}}
    t^1(e_2)(\ed \tb N\wedge\theta^2\wedge\theta^3)
    +t^2(e_2)(\ed \tb N\wedge\theta^3\wedge\theta^1)
    +t^3(e_2)(\ed \tb N\wedge\theta^1\wedge\theta^2),\\
    \ms D_{\tb u^3}I^h&=\int_{\mathcal{T}_{\mathcal{B}}}
    t^1(e_3)(\ed \tb N\wedge\theta^2\wedge\theta^3)
    +t^2(e_3)(\ed \tb N\wedge\theta^3\wedge\theta^1)
    +t^3(e_3)(\ed \tb N\wedge\theta^1\wedge\theta^2).
\end{align}
The tangent operator associated with the discrete HW functional may be formally
written as,
\begin{equation}
    \mathcal{K}=
    \begin{bmatrix}
        \ms D_{\bs \theta^i}\ms D_{\bs \theta^j}I^h & \ms D_{\bs \theta^i}\ms
        D_{\tb t^j}I^h &\ms D_{\bs \theta^i}\ms D_{\tb u ^j}I^h\\
         &\ms D_{\tb t^i}\ms D_{\tb t^j}I^h & \ms D_{\tb t^i}\ms D_{ \tb u^j}I^h\\ 
        \text{Symm}& & \ms D_{\tb u^i}\ms D_{\tb u^j}I^h
    \end{bmatrix};\quad i,j=1,...,3.
    \label{eq:tangent}
\end{equation}
The components of the tangent operator may thus be computed as,
\begin{align}
    \ms D_{\bs \theta^1 \bs \theta^1}I^h&=\int_{\mathcal{T}_{\mathcal B}}\ms D_{\bs
    \theta^1 \bs \theta^1} W \ed V;\quad
    \ms D_{\bs \theta^2 \bs \theta^2}I^h=\int_{\mathcal{T}_{\mathcal B}}\ms D_{\bs
    \theta^2 \bs \theta^2} W \ed V;\quad
    \ms D_{\bs \theta^3 \bs \theta^3}I^h=\int_{\mathcal{T}_{\mathcal B}}\ms D_{\bs
    \theta^3 \bs \theta^3} W \ed V\\
    \nonumber
    \ms D_{\bs \theta^1 \bs \theta^2}I^h&=\int_{\mathcal{T}_{\mathcal B}}\ms D_{\bs
    \theta^1 \bs \theta^2} W \ed V
    +(t^1(e_1)+t^2(e_2)+t^3(e_3))(\bs \phi \wedge \bs \phi\wedge \theta^3)\\
    &-t^3(e_1)(\bs \phi\wedge\bs \phi\wedge \ed \varphi^1)
    -t^3(e_2)(\bs \phi\wedge\bs \phi\wedge \ed \varphi^2)
    -t^3(e_3)(\bs \phi\wedge\bs \phi\wedge \ed \varphi^3)\\
    \nonumber
    \ms D_{\bs \theta^1 \bs \theta^3}I^h&=\int_{\mathcal{T}_{\mathcal B}}\ms D_{\bs
    \theta^1 \bs \theta^3} W \ed V
    +(t^1(e_1)+t^2(e_2)+t^3(e_3))(\bs \phi \wedge \bs \phi\wedge \theta^2)\\
    &-t^2(e_1)(\bs \phi\wedge\bs \phi\wedge \ed \varphi^1)
    -t^2(e_2)(\bs \phi\wedge\bs \phi\wedge \ed \varphi^2)
    -t^2(e_3)(\bs \phi\wedge\bs \phi\wedge \ed \varphi^3)\\
    \nonumber
    \ms D_{\bs \theta^2 \bs \theta^3}I^h&=\int_{\mathcal{T}_{\mathcal B}}\ms D_{\bs
    \theta^2 \bs \theta^3} W \ed V
    +(t^1(e_1)+t^2(e_2)+t^3(e_3))(\bs \phi \wedge \bs \phi\wedge \theta^1)\\
    &-t^1(e_1)(\bs \phi\wedge\bs \phi\wedge \ed \varphi^1)
    -t^1(e_2)(\bs \phi\wedge\bs \phi\wedge \ed \varphi^2)
    -t^1(e_3)(\bs \phi\wedge\bs \phi\wedge \ed \varphi^3)\\
    \ms D_{\bs \theta^1 \tb t^1}I^h&=\int_{\mathcal T_{\mathcal B} } -\bs \psi
    e_1 \otimes \bs \phi \wedge \theta^2 \wedge \theta^3;
    \ms D_{\bs \theta^2 \tb t^2}I^h=\int_{\mathcal T_{\mathcal B} } -\bs \psi
    e_2 \otimes \bs \phi \wedge \theta^3 \wedge \theta^1;
    \ms D_{\bs \theta^3 \tb t^3}I^h=\int_{\mathcal T_{\mathcal B} } -\bs \psi
    e_3 \otimes \bs \phi \wedge \theta^1 \wedge \theta^2\\
    \ms D_{\tb t^i \tb t^j}I^h&=0;\quad \ms D_{\tb u^i \tb u^j}I^h=0\\
    \ms D_{\tb u^1 \tb t^1}I^h&=\int_{\mathcal T_{\mathcal B}} \bs \psi e_1
    \otimes \ed \tb N \wedge \theta^2 \wedge \theta^3; \quad 
    \ms D_{\tb u^2 \tb t^1}I^h=\int_{\mathcal T_{\mathcal B}} \bs \psi e_2
    \otimes \ed \tb N \wedge \theta^2 \wedge\theta^3; \quad
    \ms D_{\tb u^3 \tb t^1}I^h=\int_{\mathcal T_{\mathcal B}} \bs \psi e_3
    \otimes \ed \tb N \wedge \theta^2 \wedge \theta^3\\
    \ms D_{\tb u^1 \tb t^2}I^h&=\int_{\mathcal T_{\mathcal B}} \bs \psi e_1
    \otimes \ed \tb N \wedge \theta^3 \wedge \theta^1; \quad 
    \ms D_{\tb u^2 \tb t^2}I^h=\int_{\mathcal T_{\mathcal B}} \bs \psi e_2
    \otimes \ed \tb N \wedge \theta^3 \wedge\theta^1; \quad
    \ms D_{\tb u^3 \tb t^2}I^h=\int_{\mathcal T_{\mathcal B}} \bs \psi e_3
    \otimes \ed \tb N \wedge \theta^3 \wedge \theta^1\\
    \ms D_{\tb u^1 \tb t^3}I^h&=\int_{\mathcal T_{\mathcal B}} \bs \psi e_1
    \otimes \ed \tb N \wedge \theta^1 \wedge \theta^2; \quad 
    \ms D_{\tb u^2 \tb t^3}I^h=\int_{\mathcal T_{\mathcal B}} \bs \psi e_2
    \otimes \ed \tb N \wedge \theta^1 \wedge\theta^2; \quad
    \ms D_{\tb u^3 \tb t^3}I^h=\int_{\mathcal T_{\mathcal B}} \bs \psi e_3
    \otimes \ed \tb N \wedge \theta^1 \wedge \theta^2
\end{align}
\begin{align}
    \ms D_{\bs \theta^1 \tb t^2 }I^h&=\int_{\mathcal T_{\mathcal B} }
    -\bs \psi e_2\otimes(\bs \phi \wedge \theta^2 \wedge \theta^3)
    -\bs \psi e_1\otimes(\bs \phi \wedge \theta^3 \wedge \ed \varphi^1)
    -\bs \psi e_2\otimes(\bs \phi \wedge \theta^3 \wedge \ed \varphi^2)
    -\bs \psi e_3\otimes(\bs \phi \wedge \theta^3 \wedge \ed \varphi^3)\\
    \ms D_{\bs \theta^1 \tb t^3 }I^h&=\int_{\mathcal T_{\mathcal B} }
    -\bs \psi e_3\otimes(\bs \phi \wedge \theta^2 \wedge \theta^3)
    +\bs \psi e_1\otimes(\bs \phi \wedge \theta^2 \wedge \ed \varphi^1)
    +\bs \psi e_2\otimes(\bs \phi \wedge \theta^2 \wedge \ed \varphi^2)
    +\bs \psi e_3\otimes(\bs \phi \wedge \theta^2 \wedge \ed \varphi^3)\\
    \ms D_{\bs \theta^2 \tb t^1 }I^h&=\int_{\mathcal T_{\mathcal B} }
    +\bs \psi e_1\otimes(\bs \phi \wedge \theta^1 \wedge \theta^3)
    +\bs \psi e_1\otimes(\bs \phi \wedge \theta^3 \wedge \ed \varphi^1)
    +\bs \psi e_2\otimes(\bs \phi \wedge \theta^3 \wedge \ed \varphi^2)
    +\bs \psi e_3\otimes(\bs \phi \wedge \theta^3 \wedge \ed \varphi^3)\\
    \ms D_{\bs \theta^2 \tb t^3 }I^h&=\int_{\mathcal T_{\mathcal B} }
    -\bs \psi e_3\otimes(\bs \phi \wedge \theta^3 \wedge \theta^1)
    -\bs \psi e_1\otimes(\bs \phi \wedge \theta^1 \wedge \ed \varphi^1)
    -\bs \psi e_2\otimes(\bs \phi \wedge \theta^1 \wedge \ed \varphi^2)
    -\bs \psi e_3\otimes(\bs \phi \wedge \theta^1 \wedge \ed \varphi^3)\\
    \ms D_{\bs \theta^3 \tb t^1 }I^h&=\int_{\mathcal T_{\mathcal B} }
    -\bs \psi e_1\otimes(\bs \phi \wedge \theta^1 \wedge \theta^2)
    -\bs \psi e_1\otimes(\bs \phi \wedge \theta^2 \wedge \ed \varphi^1)
    -\bs \psi e_2\otimes(\bs \phi \wedge \theta^2 \wedge \ed \varphi^2)
    -\bs \psi e_3\otimes(\bs \phi \wedge \theta^2 \wedge \ed \varphi^3)\\
    \ms D_{\bs \theta^3 \tb t^2 }I^h&=\int_{\mathcal T_{\mathcal B} }
    \bs \psi e_2\otimes(\bs \phi \wedge \theta^1 \wedge \theta^2)
    +\bs \psi e_1\otimes(\bs \phi \wedge \theta^1 \wedge \ed \varphi^1)
    +\bs \psi e_2\otimes(\bs \phi \wedge \theta^1 \wedge \ed \varphi^2)
    +\bs \psi e_3\otimes(\bs \phi \wedge \theta^1 \wedge \ed \varphi^3)\\
    \ms D_{\bs \theta^1 \bs u^1}I^h&=\int_{\mathcal T_{\mathcal B} }
    -t^2(e_1)(\ed \tb N \wedge \bs \phi \wedge \theta^3)
    +t^3(e_1)(\ed \tb N \wedge \bs \phi \wedge \theta^2)\\
    \ms D_{\bs \theta^1 \bs u^2}I^h&=\int_{\mathcal T_{\mathcal B} }
    -t^2(e_2)(\ed \tb N \wedge \bs \phi \wedge \theta^3)
    +t^3(e_2)(\ed \tb N \wedge \bs \phi \wedge \theta^2)\\
    \ms D_{\bs \theta^1 \bs u^3}I^h&=\int_{\mathcal T_{\mathcal B} }
    -t^2(e_3)(\ed \tb N \wedge \bs \phi \wedge \theta^3)
    +t^3(e_3)(\ed \tb N \wedge \bs \phi \wedge \theta^2)\\
    \ms D_{\bs \theta^2 \bs u^1}I^h&=\int_{\mathcal T_{\mathcal B} }
    +t^1(e_1)(\ed \tb N \wedge \bs \phi \wedge \theta^3)
    -t^3(e_1)(\ed \tb N \wedge \bs \phi \wedge \theta^1)\\
    \ms D_{\bs \theta^2 \bs u^2}I^h&=\int_{\mathcal T_{\mathcal B} }
    +t^1(e_2)(\ed \tb N \wedge \bs \phi \wedge \theta^3)
    -t^3(e_2)(\ed \tb N \wedge \bs \phi \wedge \theta^1)\\
    \ms D_{\bs \theta^2 \bs u^3}I^h&=\int_{\mathcal T_{\mathcal B} }
    +t^1(e_3)(\ed \tb N \wedge \bs \phi \wedge \theta^3)
    -t^3(e_3)(\ed \tb N \wedge \bs \phi \wedge \theta^1)\\
    \ms D_{\bs \theta^3 \bs u^1}I^h&=\int_{\mathcal T_{\mathcal B} }
    +t^2(e_1)(\ed \tb N \wedge \bs \phi \wedge \theta^1)
    -t^1(e_1)(\ed \tb N \wedge \bs \phi \wedge \theta^2)\\
    \ms D_{\bs \theta^3 \bs u^2}I^h&=\int_{\mathcal T_{\mathcal B} }
    +t^2(e_2)(\ed \tb N \wedge \bs \phi \wedge \theta^1)
    -t^1(e_2)(\ed \tb N \wedge \bs \phi \wedge \theta^2)\\
    \ms D_{\bs \theta^3 \bs u^3}I^h&=\int_{\mathcal T_{\mathcal B} }
    +t^2(e_3)(\ed \tb N \wedge \bs \phi \wedge \theta^1)
    -t^1(e_3)(\ed \tb N \wedge \bs \phi \wedge \theta^2)
\end{align}

\section{Numerical study}
We now demonstrate the working of the present finite element method through a few 
benchmark problems in non-linear elasticity. These test
problems are carefully designed so as to bring out any pathological behaviour
like locking and checker boarding.  There are only a few finite element
techniques in the literature that can successfully predict the deformation of
all the test problems presented in this section. To measure the convergence of
the deformation 1-forms, we use the functional norm given by,
$\int_{\mathcal{B}} \sum_{i=1}^{3}\norm{\theta^i}^2\ed V$. Similarly the
convergence of the first Piola stress is measured using
$\int_{\mathcal{B}}\norm{P}^2 \ed V$. Here $P$ is understood as a two tensor by
applying a Hodge star on the area leg.

\subsection{Cook's membrane}
\begin{figure}
    \centering
    \includegraphics[scale=0.5]{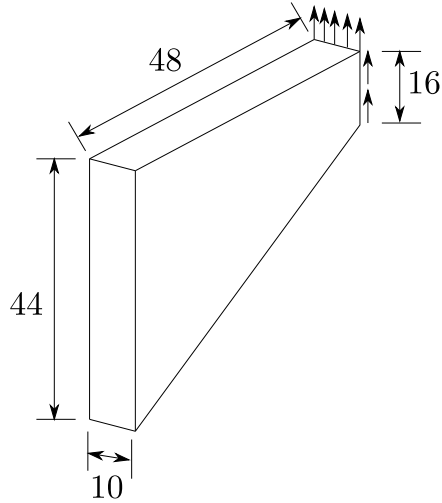}
    \caption{Dimensions of the beam, loading arrangement and boundary
    conditions for Cook's membrane problem}
    \label{fig:cookBvP}
\end{figure}

\begin{figure}
	 \centering
     \begin{subfigure}[b]{0.28\textwidth}
        \centering
        \includegraphics[scale=0.20]{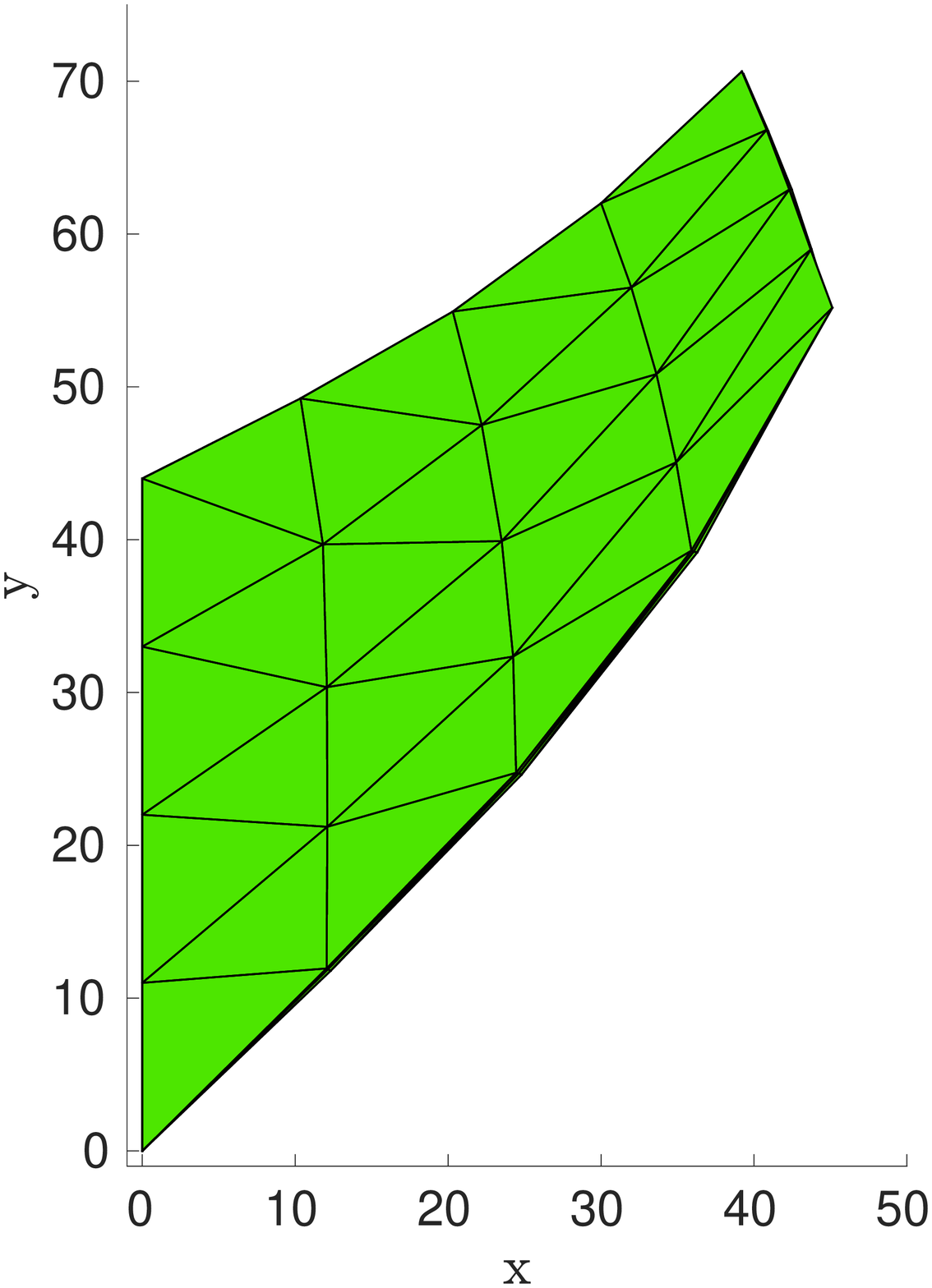}
        \caption{Mesh-1, $n=4$}
     \end{subfigure}
     \begin{subfigure}[b]{0.28\textwidth}
        \centering
        \includegraphics[scale=0.20]{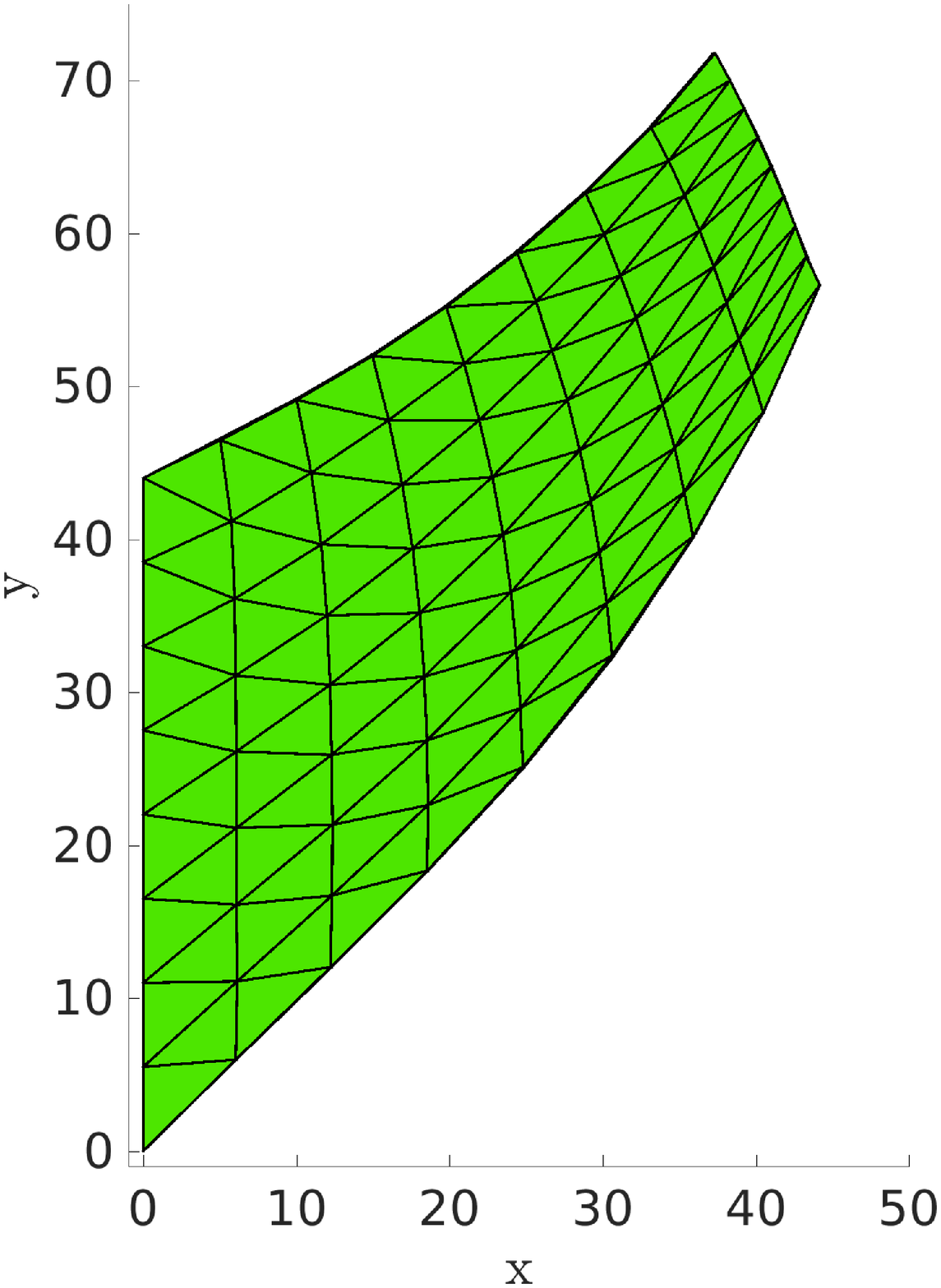}
        \caption{Mesh-2, $n=8$}
     \end{subfigure}
     \begin{subfigure}[b]{0.28\textwidth}
        \centering
        \includegraphics[scale=0.20]{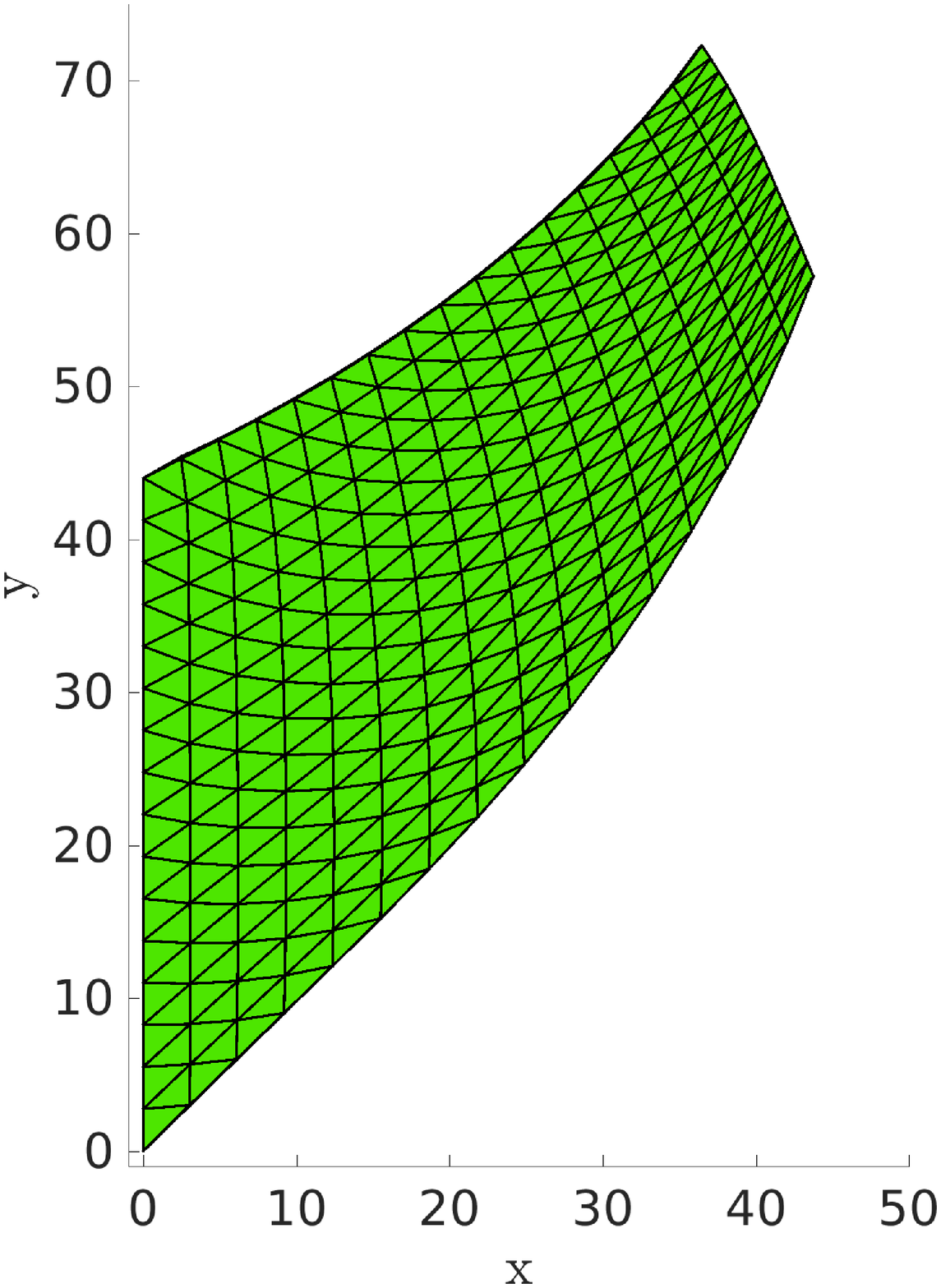}
        \caption{Mesh-3, $n=16$}
     \end{subfigure}
     \caption{Deformation of Cook's membrane predicted by our FE simulation}
     \label{fig:cookDefConfig}
\end{figure}

Cook's membrane problem is a standard benchmark to test the performance an FE
formulation against shear locking. It is now well understood that
standard displacement based FE formulations suffer from severe shear locking.
Reduced integration and enhanced strain techniques are commonly employed
to alleviate this locking, even though these techniques work well
under small deformation linear elastic assumptions. Their performance deteriorates
in large deformation regimes \cite{auricchio2013}. We test the performance of
our present FE formulation by applying it to Cook's membrane problem.  Our
objective here is to study the convergence of the finite element solution for a
non-uniform cantilever beam subjected to end shear. The dimensions of the beam,
loading arrangement and boundary conditions are is shown in Fig.
\ref{fig:cookBvP}.  A maximum load of 100 $N/mm^2$ is applied on the loading
face, which is reached in 10 loading increments.
A Mooney-Rivlin material model with the following stored energy density function
is used to describe the stress-strain relation,
\begin{equation}
    W=a(I_1-3)+b(I_2-3)+\frac{c}{2}(J-1)^2-d\log(J).
\end{equation}
In the numerical simulation, the material parameters $a$, $b$, and $c$ were set
as 126, 252, 81661, while the parameter $d$ is given by the relation
$d=2a+4b$. These material parameters correspond to a nearly incompressible
isotropic elastic material with elastic modulus  2261 and Poisson's ratio 0.4955
at the reference configuration. The contribution of the stored energy function
to the residue vector and tangent operator is thus given by,
\begin{figure}
	\centering
    \includegraphics[scale=0.45]{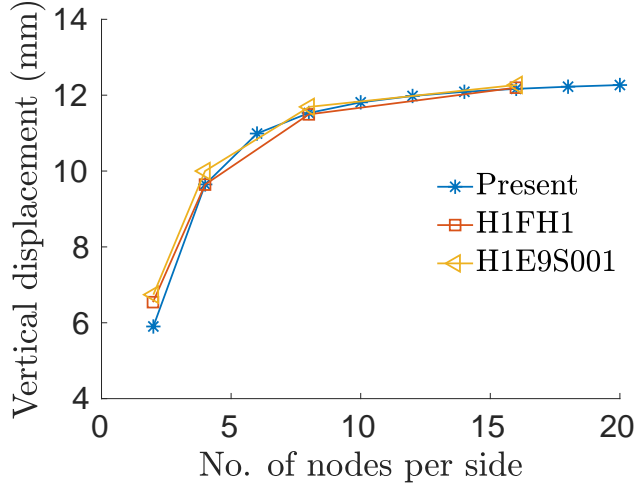}
    \caption{Convergence of tip displacement for Cook's membrane problem;
    data for finite elements H1FH1 and H1E9S001 are from Pfefferkorn and Betsch\cite{pfefferkorn2020}}
    \label{fig:cookDispConv}
\end{figure}
\begin{align}
    \ms{D}_{\bs{\theta}^i}W^h&=a\ms{D}_{\bs{\theta}^i}I_1^h+b\ms{D}_{\bs{\theta}^i}J^h+c(J^h-1)\ms{D}_{\bs{\theta}^i}J^h-\frac{d}{J^h}\ms{D}_{\bs{\theta}^i}J^h,\\
   \nonumber
    \ms{D}_{\bs{\theta}^j}\ms{D}_{\bs{\theta}^i}W^h&=a\ms{D}_{\bs{\theta}^j}\ms{D}_{\bs{\theta}^i}I_1^h
    +b\ms{D}_{\bs{\theta}^j}\ms{D}_{\bs{\theta}^i}I_2^h+c\ms{D}_{\bs{\theta}^j}J\otimes\ms{D}_{\bs{\theta}^i}J^h
    +c(J^h-1)\ms{D}_{\bs{\theta}^j}\ms{D}_{\bs{\theta}^i}J^h\\
    &-d\left[\frac{-1}{J^h}\ms{D}_{\bs{\theta}^j}J^h\otimes
    \ms{D}_{\bs{\theta}^i}J^h+\frac{1}{J^h}\ms{D}_{\bs{\theta}^j}\ms{D}_{\bs{\theta}^i}J^h
    \right].
\end{align}
The deformed configurations (for 100 $N/mm^2$ loading) predicted by our method 
for different finite element meshes are shown in Fig.
\ref{fig:cookDefConfig}. To study the convergence of tip displacement, we use a
regular mesh with $n+1$ nodes along both $x$ and $y$ directions. Fig.
\ref{fig:cookDispConv} shows the convergence of the tip displacement for
different levels of mesh refinement. Our choice of a regular mesh is dictated by
the mesh used to generate the convergence plots for the finite elements H1FH1
and H1E9S001, as given in Pfefferkorn and Betsch \cite{pfefferkorn2020}; these
finite elements use a cuboidal mesh to construct the FE spaces. In addition,
H1FH1 uses an enhancement for the deformation gradient and H1E9S001 uses an
enrichment for the cofactor of the deformation gradient.  It should also be
emphasised that, for solid mechanics problems, cuboidal finite elements are
preferred over simplicial finite elements due to their superior performance for
given number of nodes in the FE mesh.  Moreover, it may be seen from Fig.
\ref{fig:cookDispConv}, our simplicial finite element has comparable performance 
vis-\'a-vis the cuboidal finite elements. 
\begin{figure}
	 \centering
     \begin{subfigure}[b]{0.48\textwidth}
        \centering
        \includegraphics[scale=0.25]{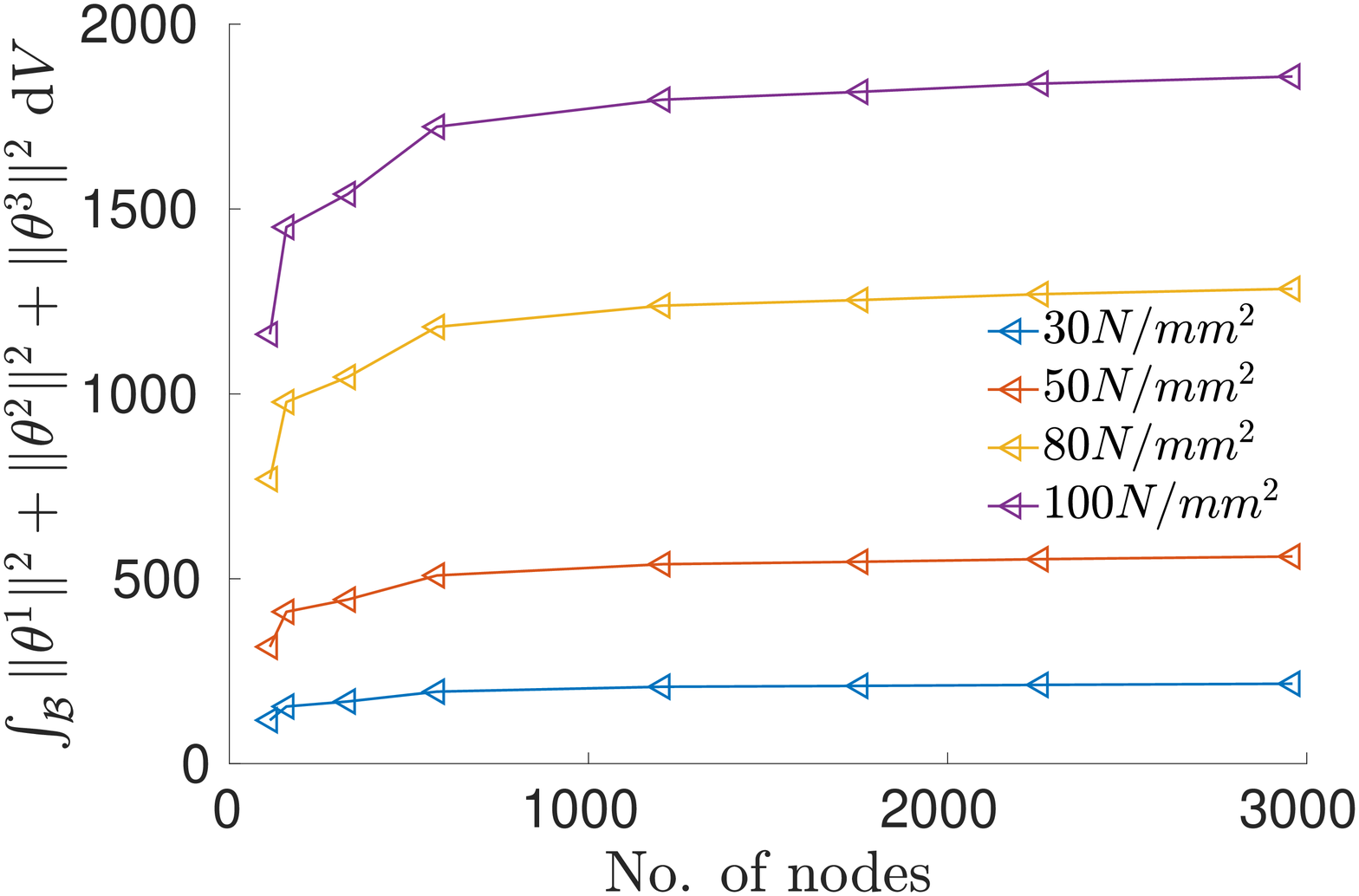}
        \caption{Convergence of deformation 1-forms}
        \label{fig:cookThetaConv}
     \end{subfigure}
     \begin{subfigure}[b]{0.48\textwidth}
        \centering
        \includegraphics[scale=0.25]{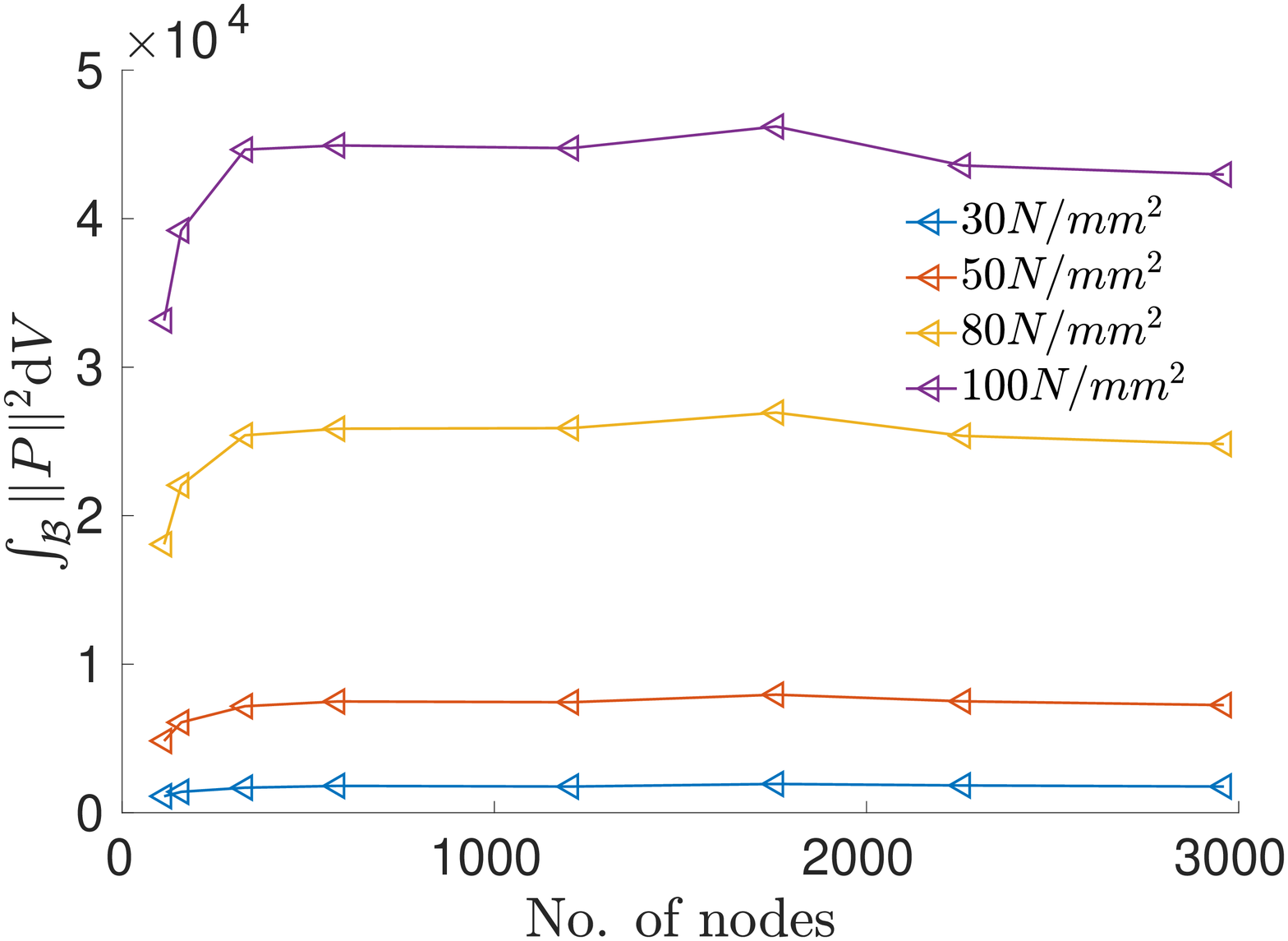}
        \caption{Convergence of first Piola stress}
        \label{fig:cookStressConv}
     \end{subfigure}
     \caption{Convergence of deformation 1-forms and first Piola stress for Cook's
     membrane problem}
     \label{fig:cookConvField}
\end{figure}
The convergence of deformation 1-forms and first Piola stress are shown in Fig.
\ref{fig:cookConvField}. To generate these convergence plots, FE meshes with
refinement near the root of the beam were used. Regular meshes with uniform
refinement used to study the convergence of displacements required a
prohibitively large computational cost to resolve the higher gradients near the root.

\subsection{Cube compression}
\begin{figure}
    \centering
    \includegraphics[scale=1.0]{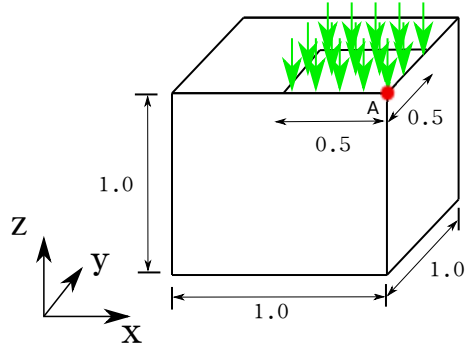}
    \caption{Quarter model of an elastic cube under compression}
    \label{fig:ccBvP}
\end{figure}

This problem is designed to test the efficiency of the FE approximation against
numerical instabilities like hour glassing \cite{reese2000,shojaei2019} under
extremely large deformation. Conventional displacement based finite elements
suffer from such instabilities when subjected to large deformation. In the
literature, enhanced  strain techniques -- again a family of numerical fixes -- 
have been suggested to alleviate such
instabilities. We now apply the proposed FE formulation to a cuboid subjected to
extremely large compression. The cuboid has a side of length 2 units. To reduce
the computational effort, we model only a quarter of the domain.  A quarter of
the domain along with the loading condition is shown in Fig.  \ref{fig:ccBvP}.
This simplification is facilitated by the symmetry of the loading and the
boundary conditions. A uniform load is applied on a square region whose side is 1
unit on the top face of the cuboid. This translates to a square of side 0.5
units for the quarter domain. The bottom face of the cuboid is constrained to
move vertically, while the loaded region is prevented to move horizontally. A
maximum load of $320 N/mm^2$ was applied in the loaded region, which was
reached in multiple load increments. In addition to these, zero displacement
boundary conditions are imposed to enforce symmetry of deformation. A
Neo-Hookean stored energy function is adopted to represent the stress-strain 
behaviour of the block. The stored energy function of the material model is given
by,
\begin{equation}
W=\frac{\mu}{2}(I_1-3)-\frac{\mu}{2}\log (J)+\frac{\kappa}{8}\log (J)^2.
\end{equation}
In the numerical simulation, the material parameters $\lambda$ and $\mu$ were
set as 400889.806 and 80.194 respectively. These material constants correspond
to an isotropic quasi-incompressible material in the reference configuration. The
first and the second derivatives of the stored energy function with respect to
the degrees of freedom associated with deformation 1-forms may be computed as,
\begin{align}
\ms{D}_{\bs{\theta}^i}W^h&=\frac{\mu}{2}\ms{D}_{\bs\theta^i}I_1^h-\frac{\mu}{2J^h}\ms{D}_{\bs\theta^i}J^h
+\frac{\kappa}{4J^h}\log J^h\ms{D}_{\bs\theta^i}J^h,\\
\ms{D}_{\bs{\theta}^i}\ms{D}_{\bs{\theta}^j}W^h&=
    \frac{\mu}{2}\ms{D}_{\bs{\theta}^i}\ms{D}_{\bs{\theta}^j}I_1^h+\left(\frac{\mu}{2(J^h)^2}
    +\frac{\kappa}{4(J^h)^2}-\frac{\kappa \log J^h}{4(J^h)^2} \right)\ms{D}_{\bs{\theta}^i}J^h \otimes\ms{D}_{\bs{\theta}^j}J^h
+\left(-\frac{\mu}{2J^h}+\frac{\kappa}{4J^h}
    \right)\ms{D}_{\bs{\theta}^i}\ms{D}_{\bs{\theta}^j}J^h.
\end{align}

\begin{figure}
     \centering
     \begin{subfigure}[b]{0.48\textwidth}
         \centering
         \includegraphics[scale=0.45]{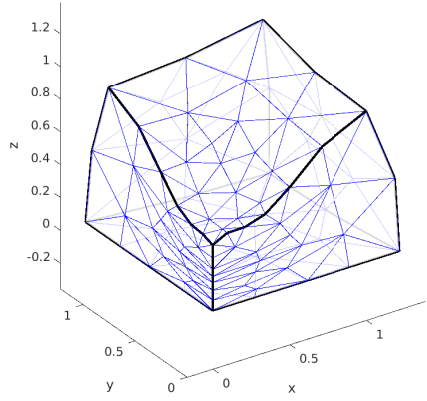}
         \caption{Mesh-1}
     \end{subfigure}
     \begin{subfigure}[b]{0.48\textwidth}
         \centering
         \includegraphics[scale=0.55]{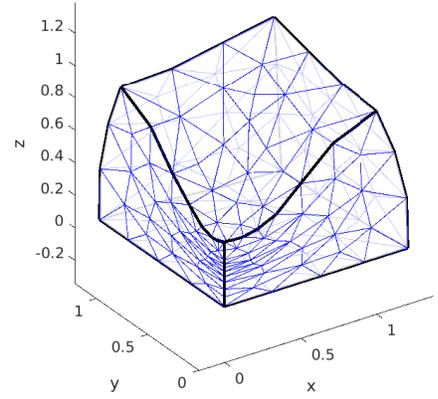}
         \caption{Mesh-2}
     \end{subfigure}
     \begin{subfigure}[b]{0.48\textwidth}
         \centering
         \includegraphics[scale=0.45]{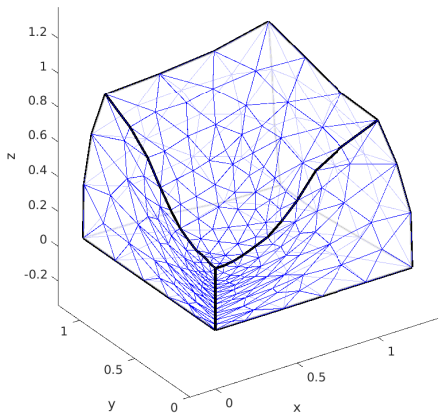}
         \caption{Mesh-3}
     \end{subfigure}
     \begin{subfigure}[b]{0.48\textwidth}
         \centering
         \includegraphics[scale=0.45]{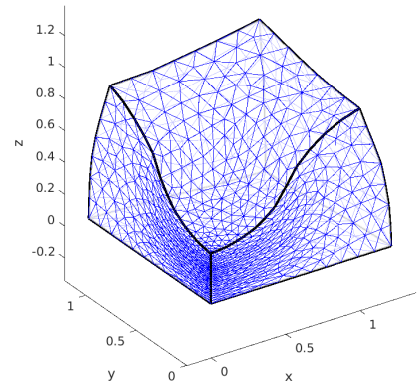}
         \caption{Mesh-4}
     \end{subfigure}
        \caption{Deformed shapes of the cube predicted by our FE approximation
        for different finite element meshes}
        \label{fig:ccDefConf}
\end{figure}

\begin{figure}
    \centering
    \includegraphics[scale=0.40]{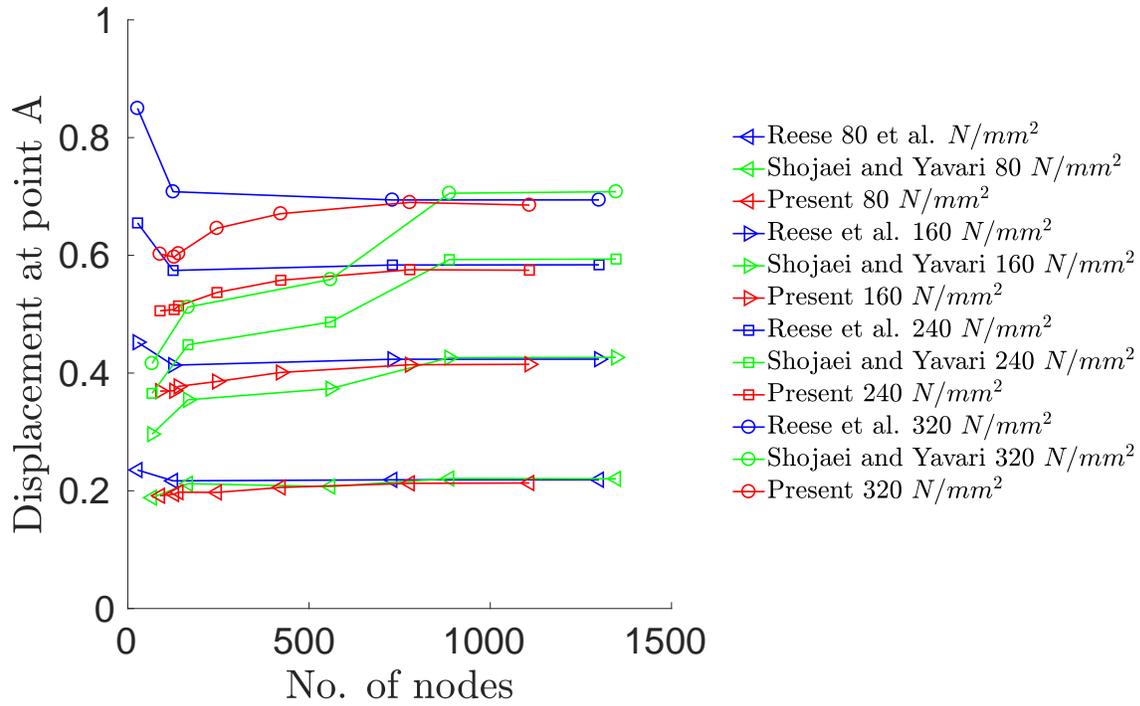}
    \caption{Convergence of displacement at the centre of the loading face}
    \label{fig:ccDispConv}
\end{figure}

\begin{figure}
		\centering
    \begin{subfigure}[b]{0.48\textwidth}
        \centering
        \includegraphics[scale=0.25]{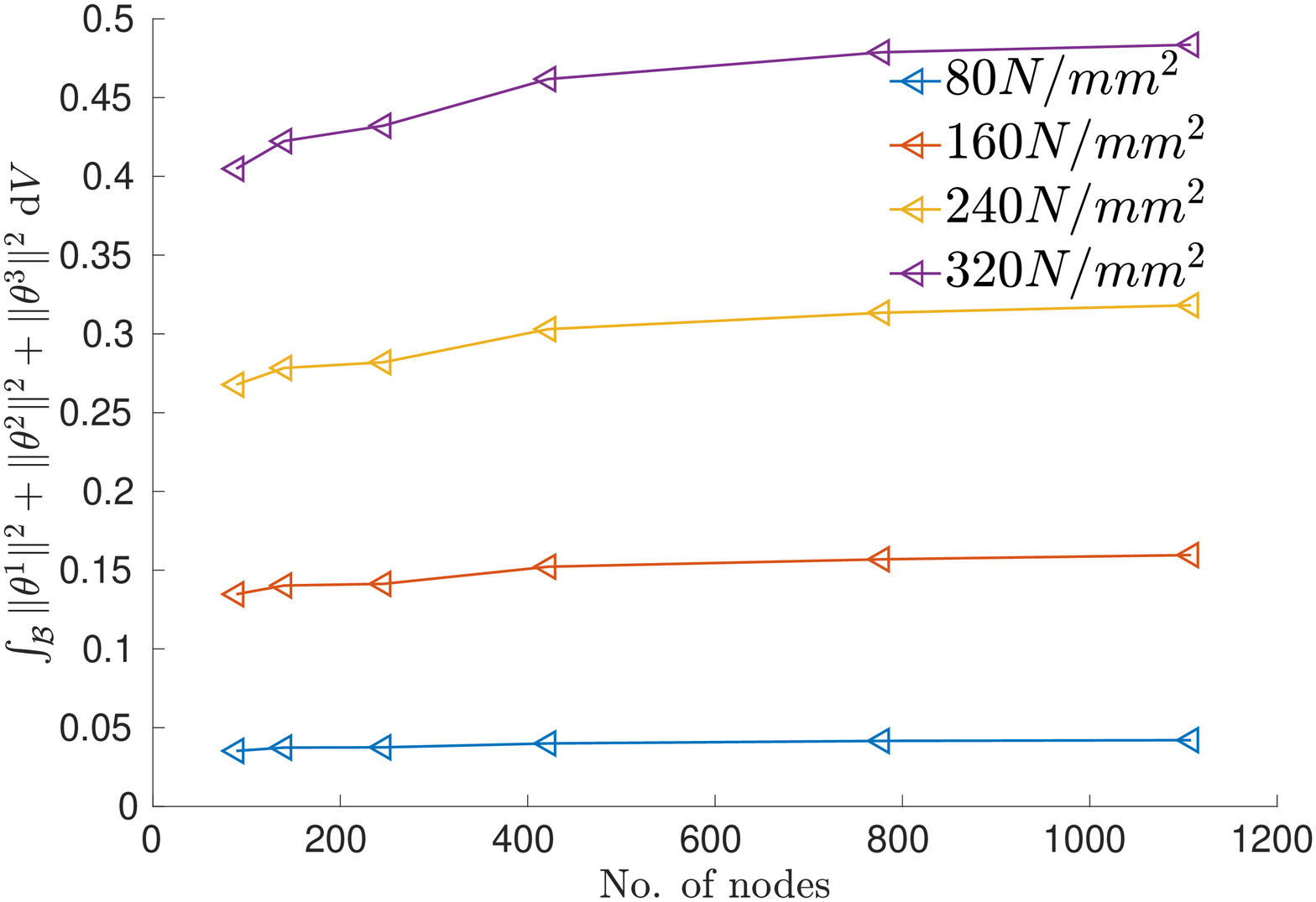}
        \caption{Convergence of deformation 1-forms}
        \label{fig:ccThetaConv}
    \end{subfigure}
    \begin{subfigure}[b]{0.48\textwidth}
        \centering
        \includegraphics[scale=0.25]{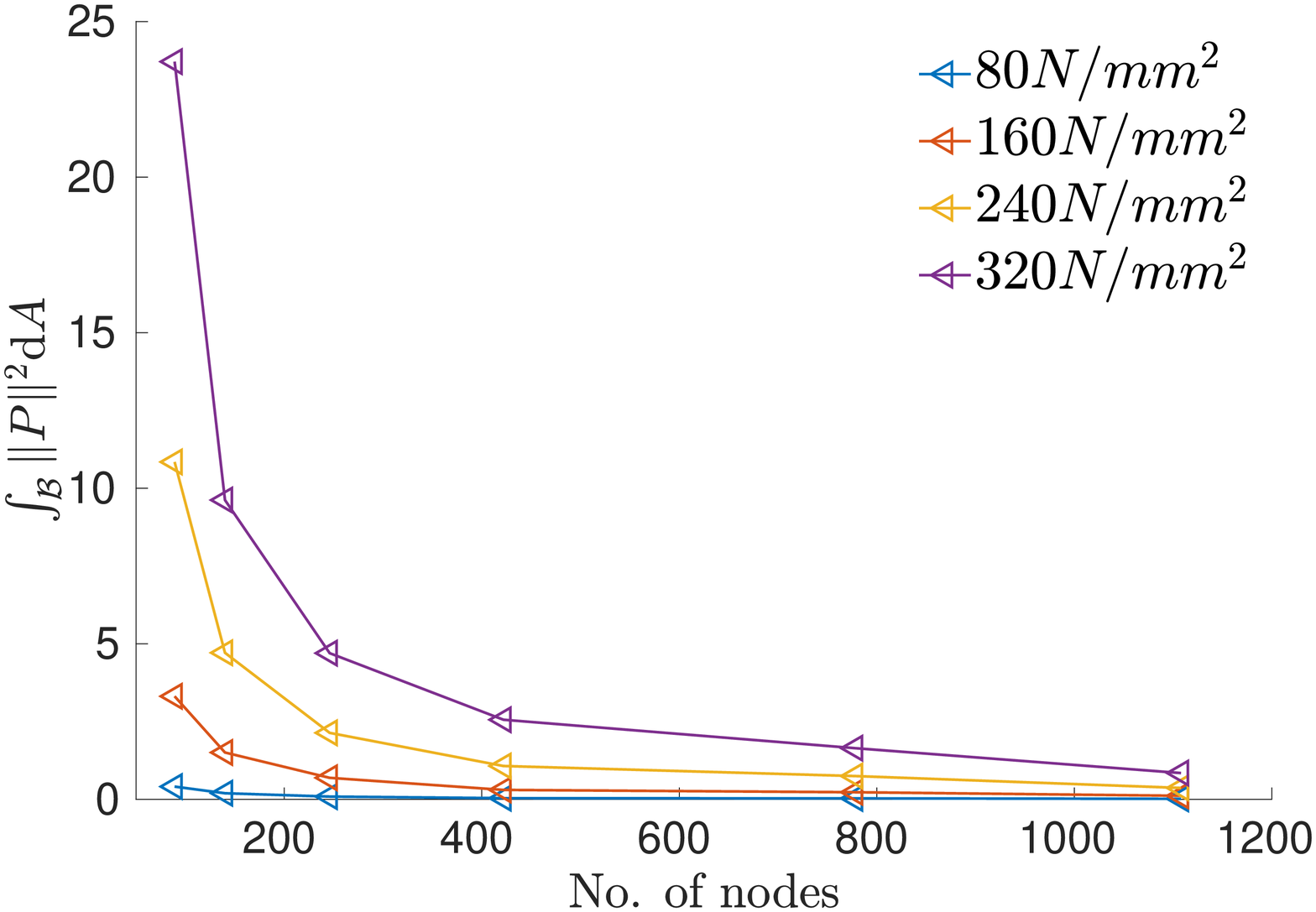}
        \caption{Convergence of Stress}
        \label{fig:ccStressConv}
    \end{subfigure}
    \caption{Convergence of deformation 1-form and first Piola stress for the
    cube compression problem}
    \label{fig:ccConvField}
\end{figure}

The deformation predicted by our FE formulation for different simplicial meshes is
shown in Fig. \ref{fig:ccDefConf}. For a quantitative comparison, we show the
convergence of displacement in the $z$ direction at point A in Fig.
\ref{fig:ccDispConv}, where the prediction of our FE formulation is plotted along
with those by Reese \textit{et al.} \cite{reese2000} and Shojaei and
Yavari \cite{shojaei2019}. The solution technique adopted by Shojaei and Yavari
is a mixed finite element technique with an additional stabilization term which
vanishes at the equilibrium point. Reese \textit{et al.}, on the other hand, have used an
enhanced strain based stabilisation technique to arrive at their prediction.
From the displacement convergence plots, it is seen that our formulation
performs comparable to previously mentioned techniques. However, the main
advantage of out formulation is that no additional parameters (like
stabilisation parameters) other that the material constants are required in the 
numerical procedure. Convergence plots of deformation 1-forms and the first Piola
stress are given in Fig. \ref{fig:ccConvField}. A
near monotonic convergence is observed. 

\subsection{Torsion} 
\begin{figure}
    \centering
    \includegraphics[scale=0.7]{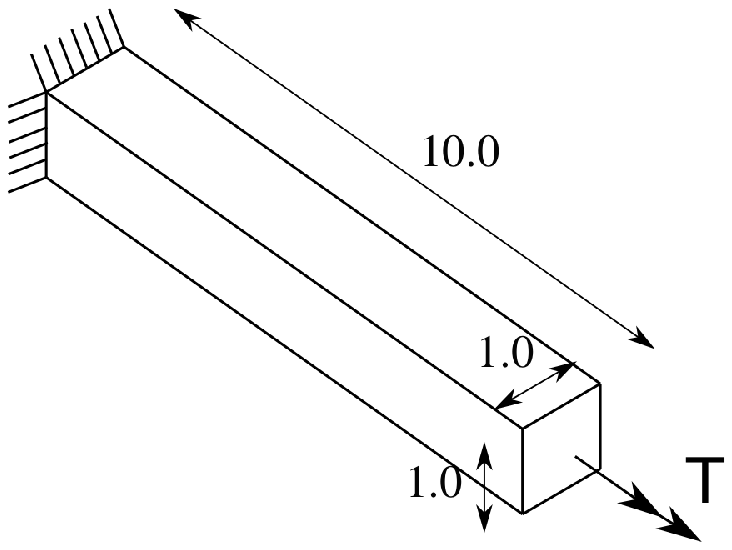}
    \caption{Dimensions, boundary conditions and loading arrangement of a square
    rod subjected to torsion}
    \label{fig:torsionBvP}
\end{figure}

\begin{figure}
     \centering
     \begin{subfigure}[b]{0.48\textwidth}
         \centering
         \includegraphics[scale=0.25]{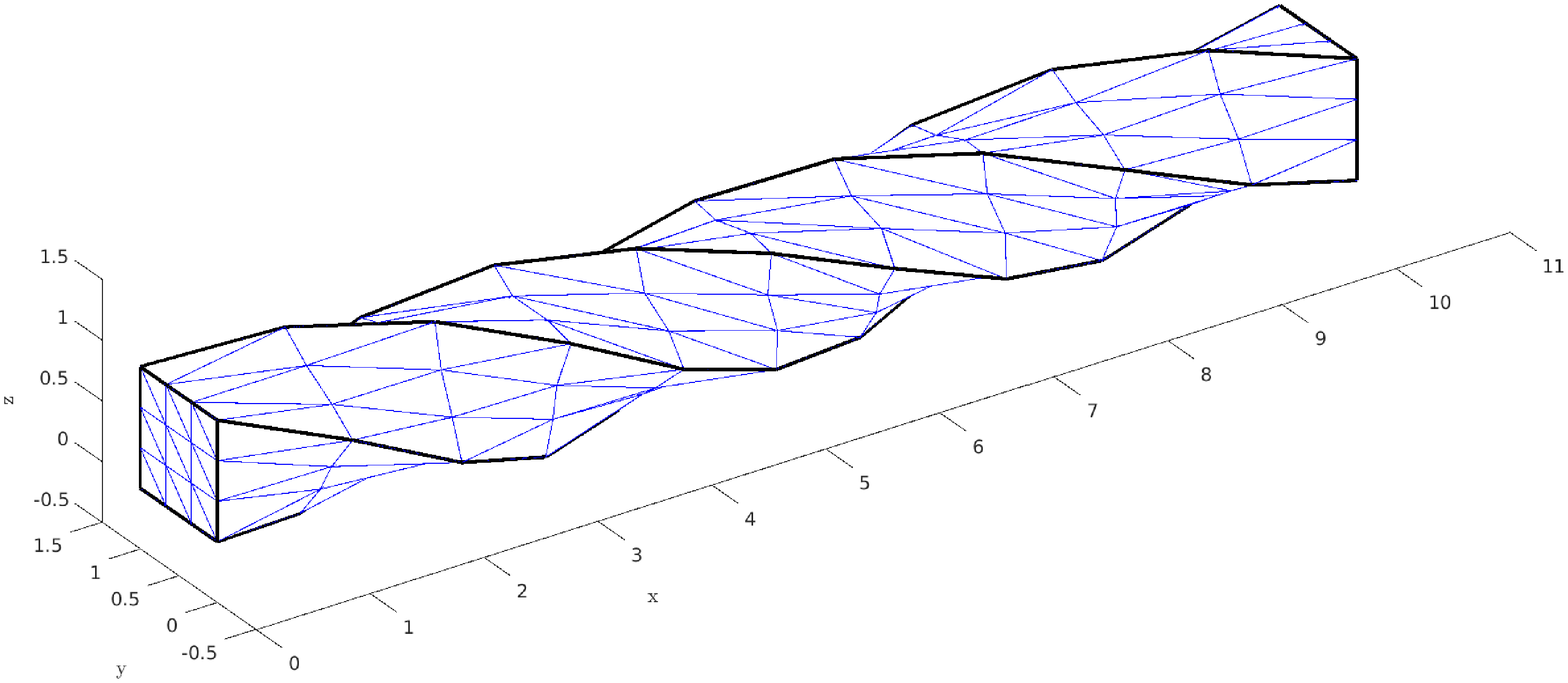}
         \caption{Mesh-1}
         \label{fig:torsion1}
     \end{subfigure}
     \begin{subfigure}[b]{0.48\textwidth}
         \centering
         \includegraphics[scale=0.30]{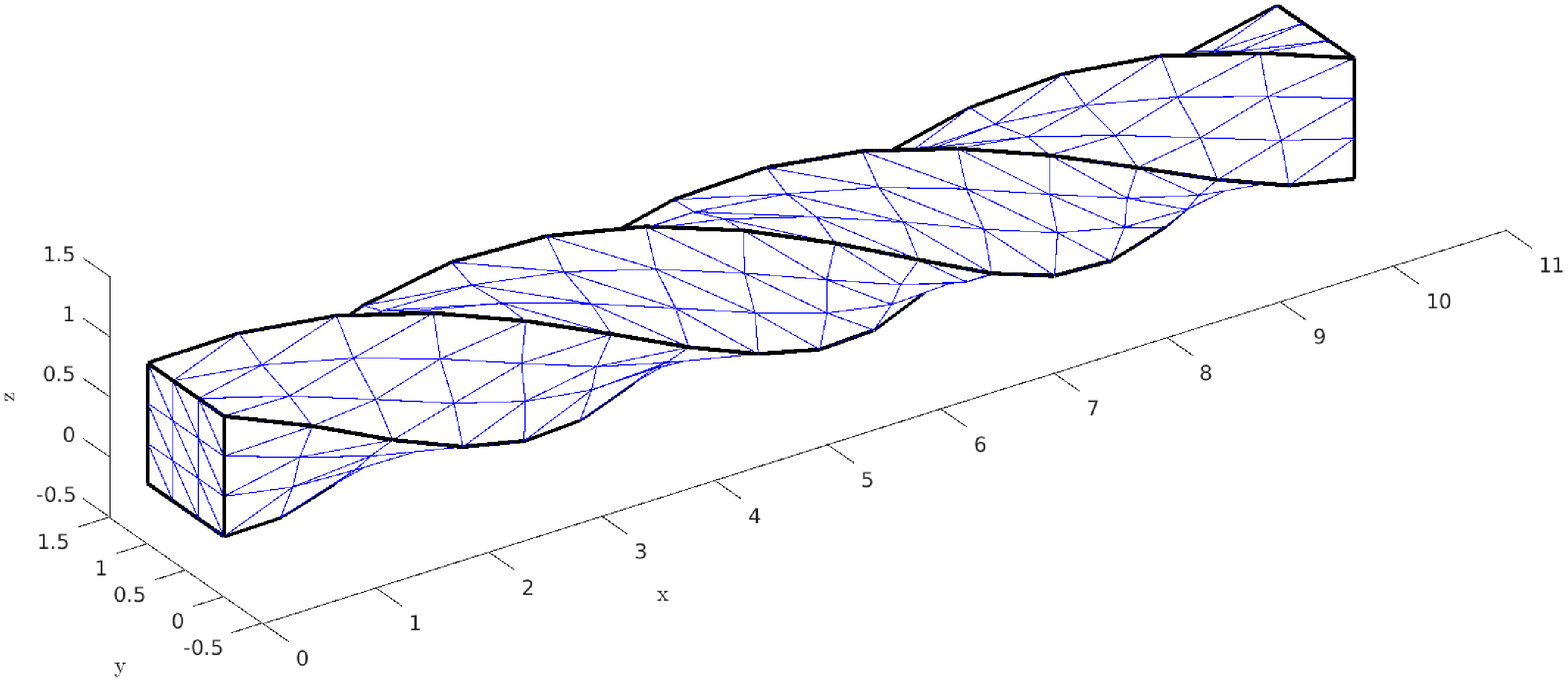}
         \caption{Mesh-2}
         \label{fig:torsion2}
     \end{subfigure}\\
     \centering
     \begin{subfigure}[b]{0.48\textwidth}
         \centering
         \includegraphics[scale=0.30]{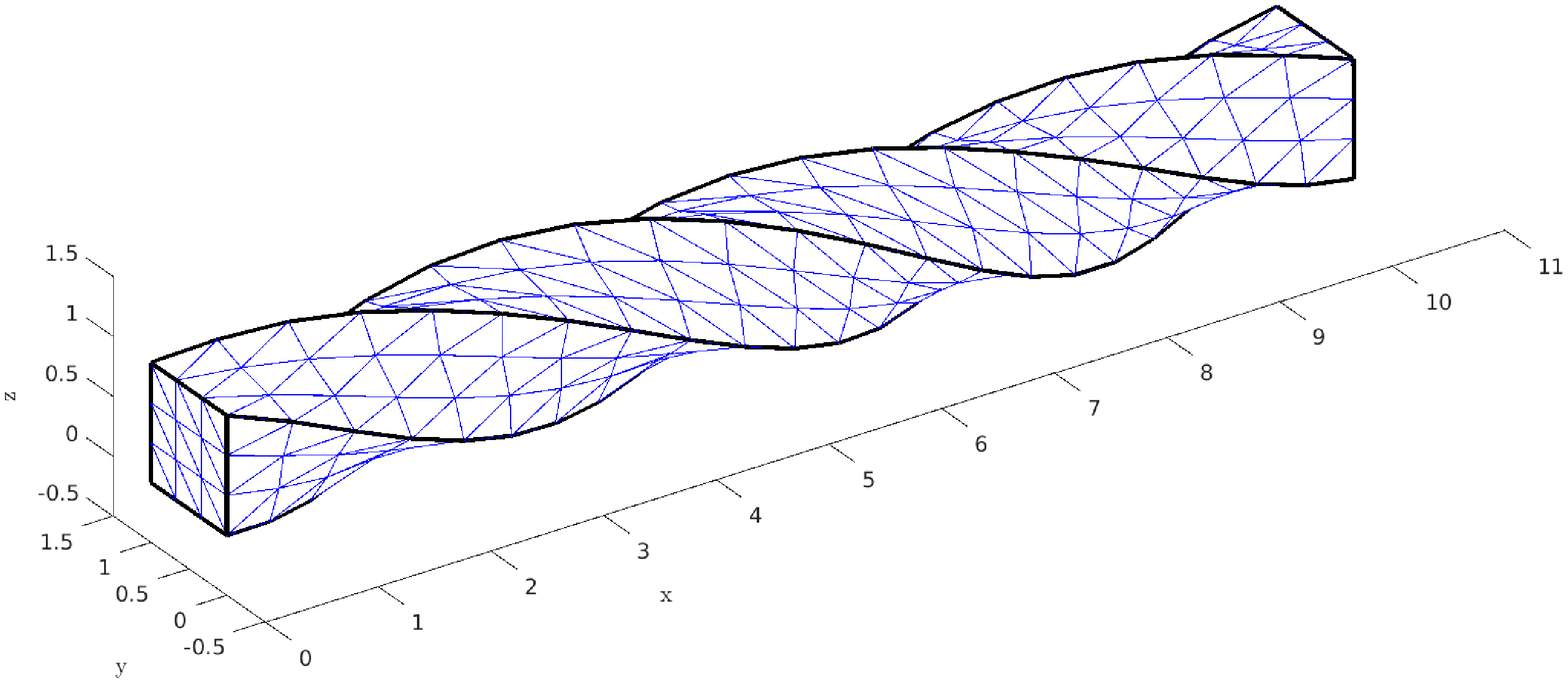}
         \caption{Mesh-3}
         \label{fig:torsion3}
     \end{subfigure}
     \begin{subfigure}[b]{0.48\textwidth}
         \centering
         \includegraphics[scale=0.30]{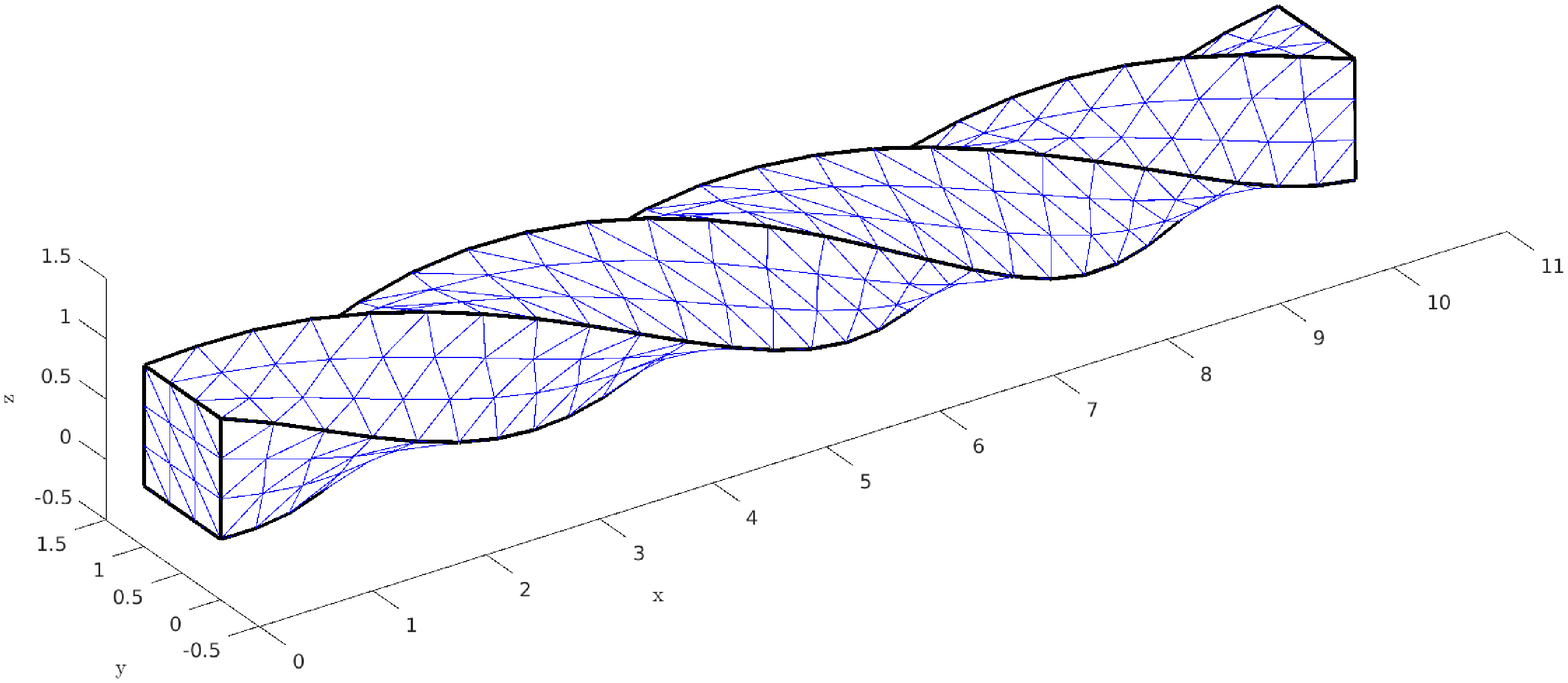}
         \caption{Mesh-4}
         \label{fig:torsion4}
     \end{subfigure}
     \caption{Deformed shapes of the cube computed using four different finite element meshes}
     \label{fig:torsionDefConfig}
\end{figure}

\begin{figure}
     \centering
     \begin{subfigure}[b]{0.49\textwidth}
        \centering
        \includegraphics[scale=0.29]{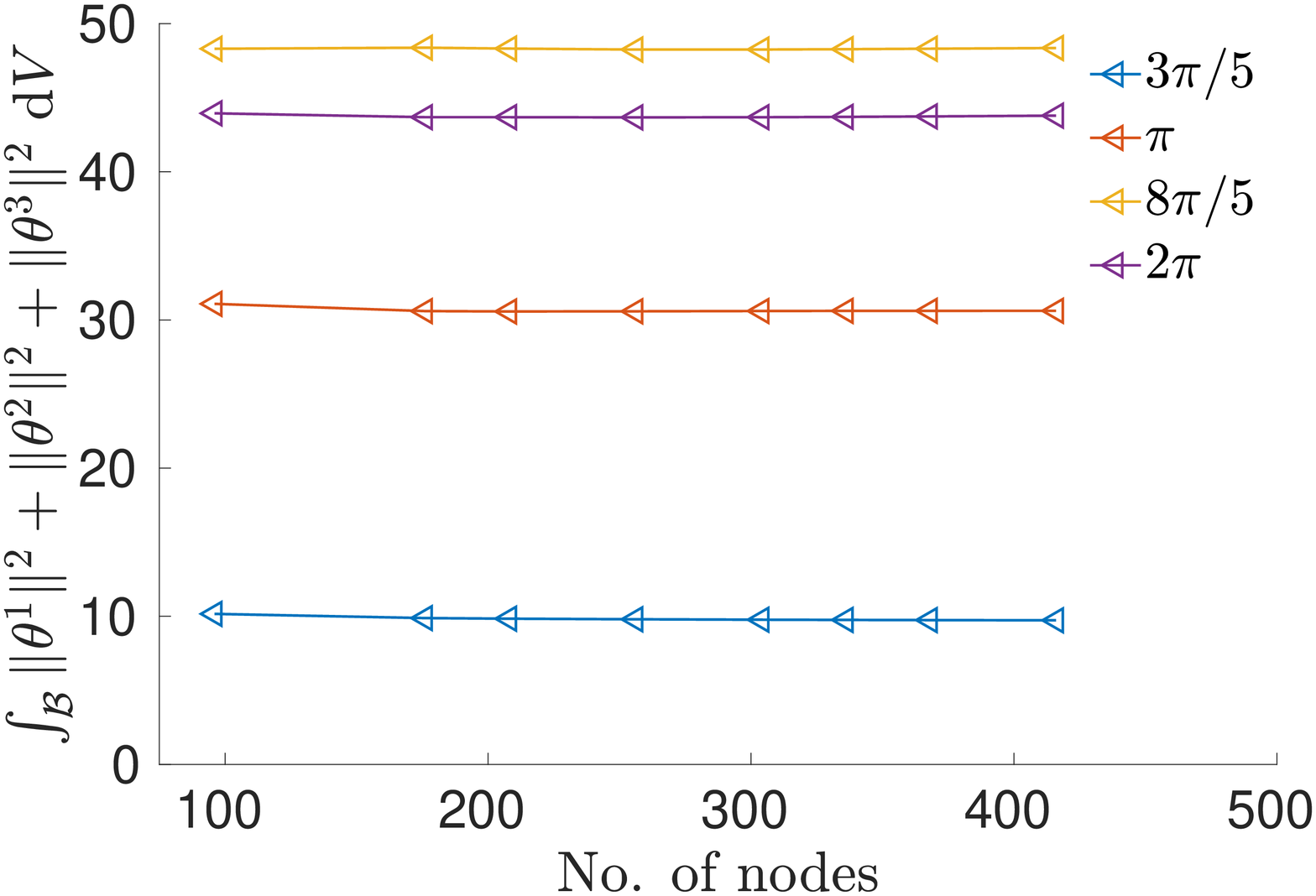}
         \caption{Convergence of deformation 1-form}
     \end{subfigure}
     \begin{subfigure}[b]{0.49\textwidth}
        \centering
        \includegraphics[scale=0.29]{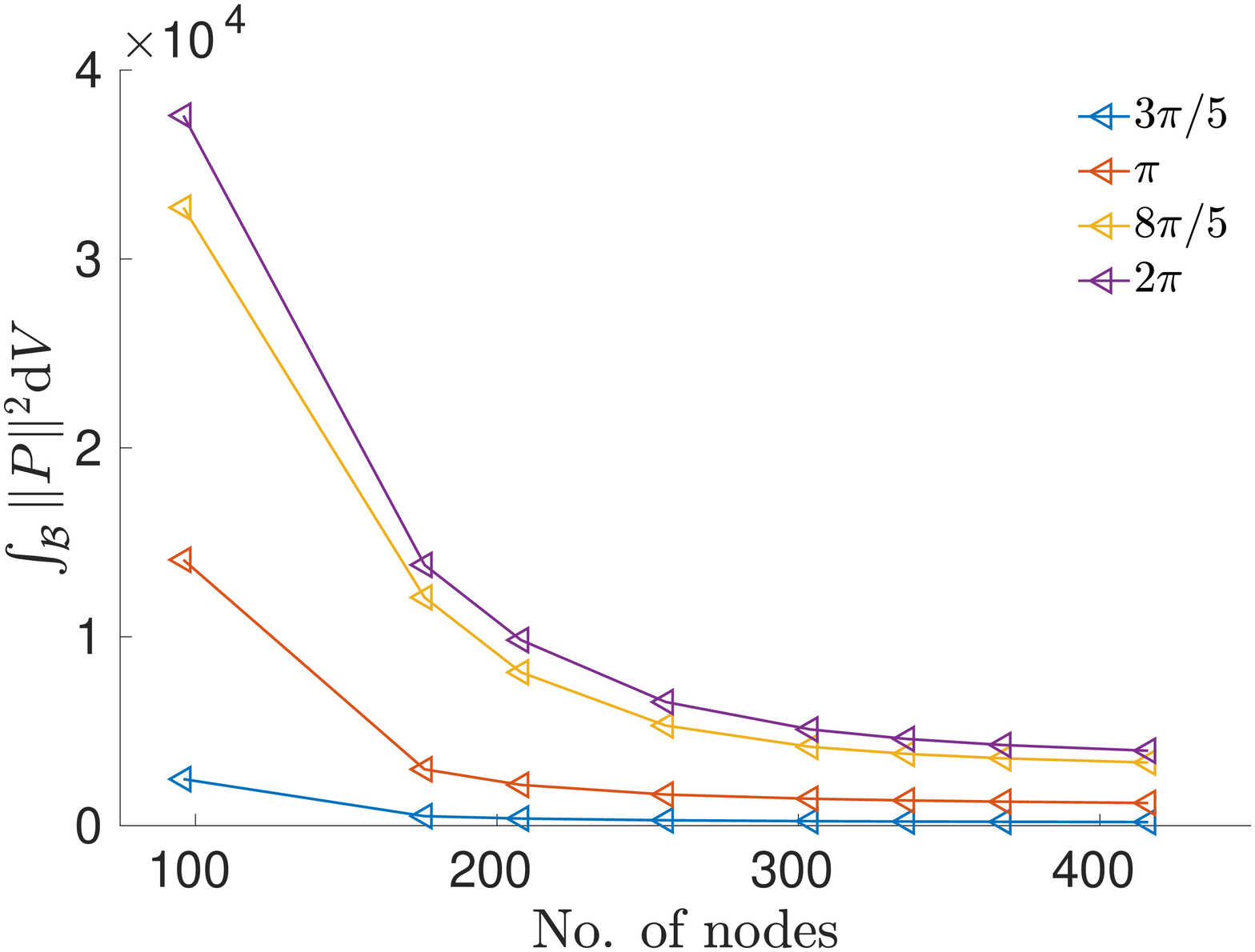}
         \caption{Convergence of traction 1-form}
     \end{subfigure}
     \caption{Convergence of deformation and traction 1-forms for the torsion
     problem}
     \label{fig:torsionConvField}
\end{figure}

We now study the effect of extreme mesh distortion during the deformation
process. A solid shaft with 1 x 1 cross section and length 10 units is held at
one end and subjected to a twist on the other. Fig \ref{fig:torsionBvP}
shows the domain of the problem along with the boundary conditions. The left end
of the shaft is constrained in all directions, while on the twisting face
all out-of-plane motions are constrained. Instead of applying a twisting moment at the
twisting face, displacements are prescribed.  In other words, the simulation
is performed by rotation control (of the twisting face) and not by a twisting
moment. A maximum of 2$\pi$ rotation is applied at the twisting face; this
rotation is reached in multiple rotation increments.  A compressible isotropic
Mooney-Rivilin type  material model is used to describe the stress strain
relation, whose stored energy function is given by,
\begin{equation}
    W=\frac{\alpha}{2}I_1^2+\frac{\beta}{2}(I_2)^2-\nu\log (J)
\end{equation}
The material constants $\alpha$ and $\beta$ were chosen as 24 and 84
respectively and $\nu=(6\alpha+12 \beta)$. The contributions of the stored
energy function to the residue vector and tangent operator are given by,
\begin{align}
\ms{D}_{\bs{\theta}^i}W^h=&
\alpha I_1^h\ms{D}_{\bs\theta^1}I_1^h+\beta I_2^h\ms{D}_{\bs\theta^i}
-\frac{\nu}{J^h}\ms{D}_{\bs{\theta}^i}J^h,\\
\nonumber
\ms{D}_{\bs{\theta}^i}\ms{D}_{\bs{\theta}^j}W^h=&\alpha(\ms{D}_{\bs\theta^i}I_1^h\otimes\ms{D}_{\bs\theta^j}I_1^h+I_1\ms{D}_{\bs\theta^i}\ms{D}_{\bs\theta^j}I_1^h)
+\beta(\ms{D}_{\bs\theta^j}I_2^h\otimes \ms{D}_{\bs\theta^j}I_2^h+I_2^h\ms{D}_{\bs{\theta}^i}\ms{D}_{\bs{\theta}^j}I_1^h)\\
    &-\nu\left(-\frac{1}{(J^h)^2}\ms{D}_{\bs\theta^i}J^h
    \otimes\ms{D}_{\bs\theta^j}J^h+\frac{1}{J^h}\ms{D}_{\bs\theta^i}\ms{D}_{\bs\theta^j}J^h
    \right).
\end{align}
The deformed configurations predicted by our FE procedure for different FE meshes
and for a twist angle of $2\pi$ are shown in Fig. \ref{fig:torsionDefConfig}. The
convergence of deformation 1-form and the first Piola  stress for different twist
angles is shown in Fig. \ref{fig:torsionConvField}.

\subsection{Shearing of split ring}
\begin{figure}
    \centering
    \includegraphics[scale=0.95]{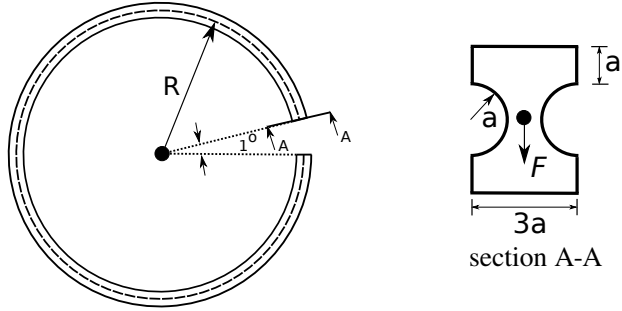}
    \caption{Geometry and loading condition for a split ring under shear}
    \label{fig:srBvP}
\end{figure}

\begin{figure}
    \centering
    \begin{subfigure}[b]{0.49\textwidth}
        \includegraphics[scale=0.25]{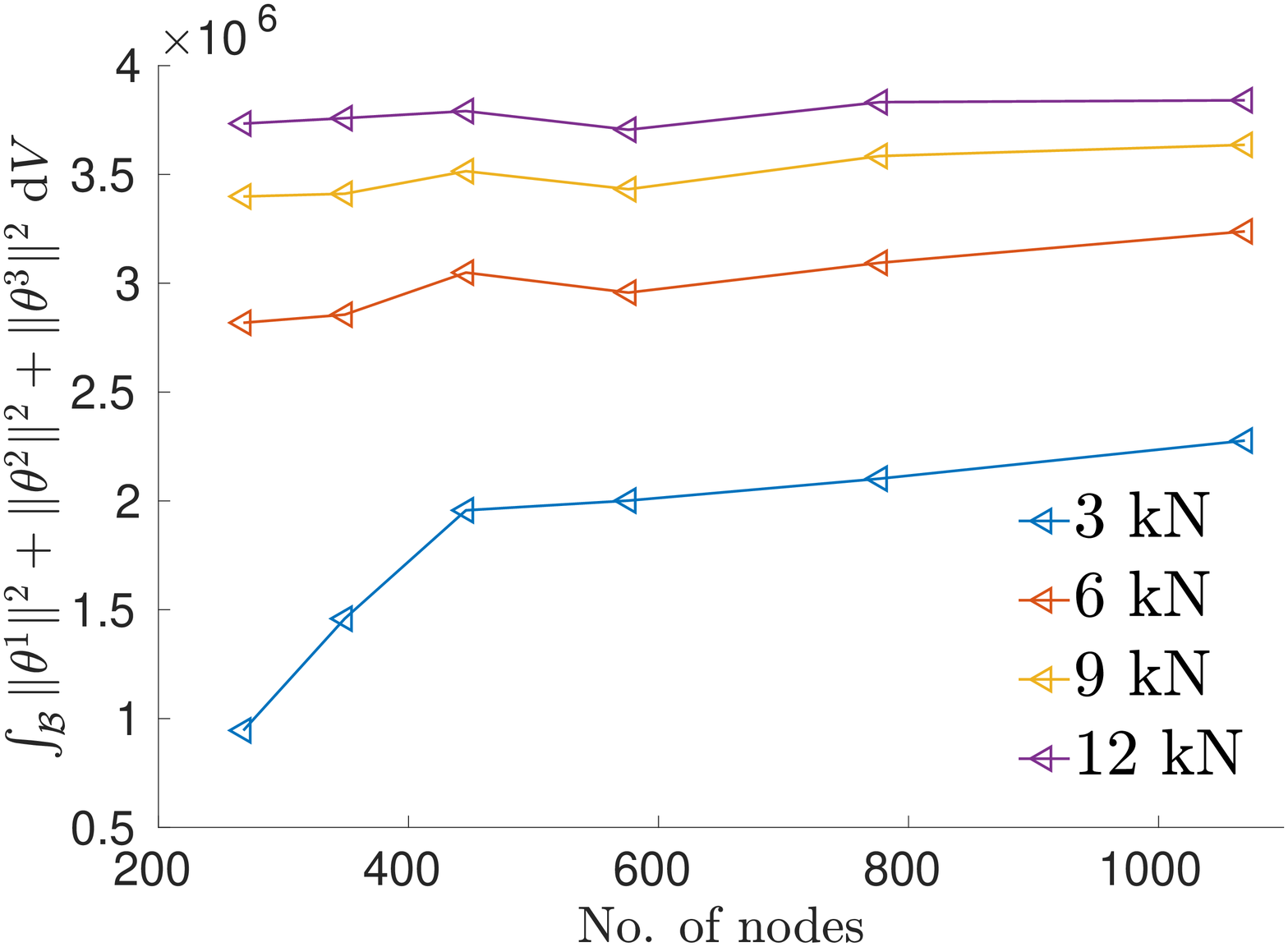}
        \caption{Convergence of deformation 1-forms}
        \label{fig:srThetaConv}
    \end{subfigure}
    \begin{subfigure}[b]{0.49\textwidth}
        \includegraphics[scale=0.25]{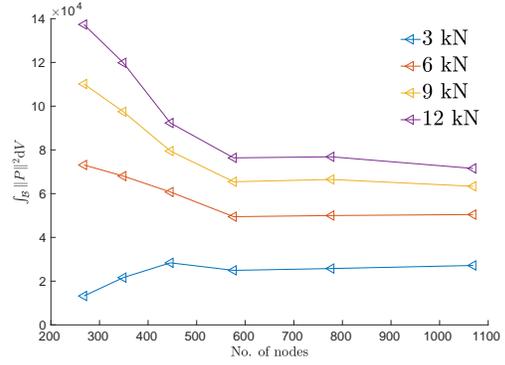}
        \caption{Convergence of traction 1-forms}
        \label{fig:srStressConv}
    \end{subfigure}
    \caption{Convergence of deformation and traction 1-forms for the split ring
    problem}
    \label{fig:srConvField}
\end{figure}

\begin{figure}
     \centering
     \begin{subfigure}[b]{0.30\textwidth}
         \centering
         \includegraphics[scale=0.38]{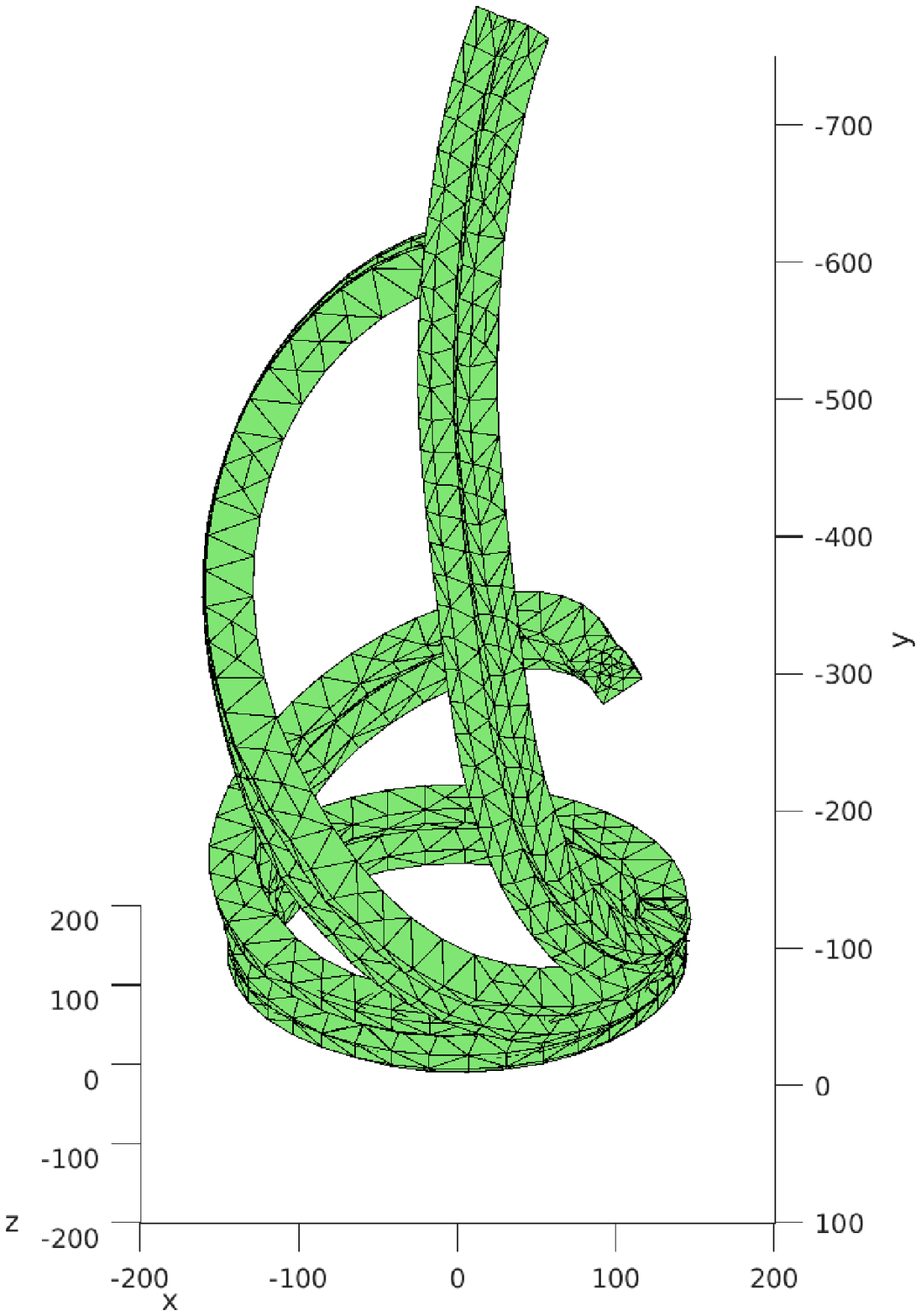}
         \caption{Mesh-2}
     \end{subfigure}
     \begin{subfigure}[b]{0.30\textwidth}
         \centering
         \includegraphics[scale=0.38]{./images/srDeformedConfigMesh2.eps}
         \caption{Mesh-3}
     \end{subfigure}
          \begin{subfigure}[b]{0.30\textwidth}
         \centering
         \includegraphics[scale=0.38]{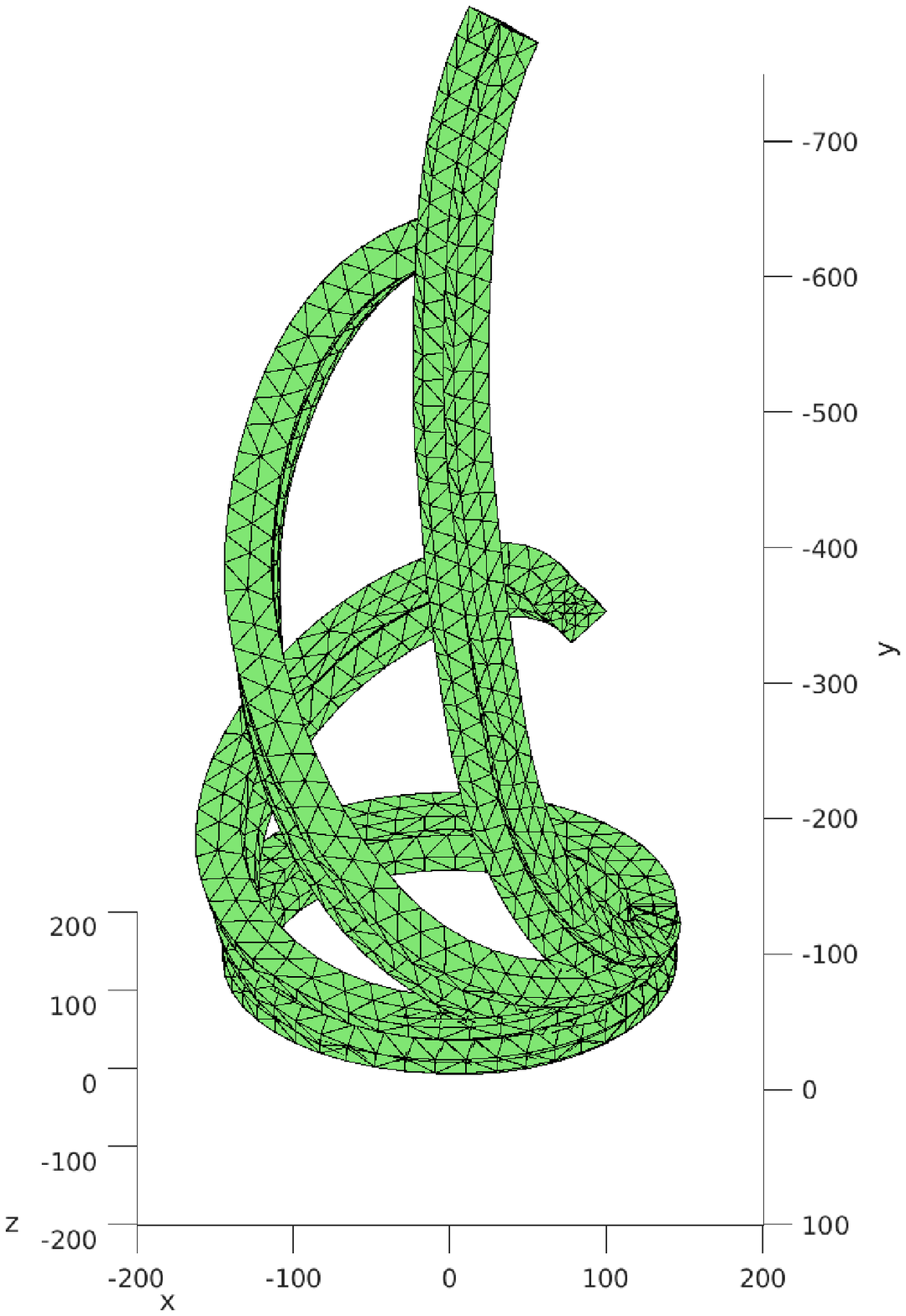}
         \caption{Mesh-3}
     \end{subfigure}
        \caption{Deformed shapes of the split ring at different load steps}
        \label{fig:srDefConfig}
\end{figure}

Finally, we assess the performance of our formulation for a problem where the domain is a little
complicated and the deformation is extremely large, even though strains are small (large
rotations in some parts of the domain). The domain is an incomplete annulus,
obtained by revolving an I shaped cross section to 359 degree about the $y$-axis
(the cross section at 0 degree revolution is contained in the $x-z$ plane). The
cross section of the annulus at zero degree revolution angle is constrained against
all motion. A shear force along the $y$-direction is applied at the centroid of
the cross section located at a revolution angle of 359 degree. The cross section
and boundary conditions are illustrated in Fig.  \ref{fig:srBvP}. Wackerfu{\ss}
and Gruttmann \cite{wackerfuss2011} used a rod formulation to study the
deformation of this problem.  For the geometry and loading conditions, it may be
observed that cross sections along the length of the ring undergo a combination
of bending, twisting, shearing and extension.  In addition, a large segment of
the ring undergoes rigid rotation with small strains.  A maximum shear force of 12
kN is applied which is reached in multiple loading increments. A Neo-Hookian type material
model is used where the stored energy function is given by,
\begin{equation}
    W=\frac{\mu}{2}(I_1-3)-\mu \ln J + \frac{\lambda}{4}(J^2-1-2\ln J)
    \label{eq:WSplitRing}
\end{equation}
The material
parameters are takes as $\mu =0.4$ and $\lambda=4\alpha$. For the assumed
stored energy function, contributions to the residue vector and tangent operator
are given by,
\begin{align}
    \ms{D}_{\bs\theta^i}W^h&=\frac{\mu}{2}\ms{D}_{\bs\theta^1}I_1-\frac{\mu}{J}\ms
    D_{\bs \theta^i}+\frac{\lambda}{4}(2J^h\ms D_{\bs
    \theta^i}J^h-\frac{2}{J}\ms D_{\bs \theta^i}J^h)\\
    \nonumber
    \ms{D}_{\bs \theta^i} \ms{D}_{\bs \theta^j}W^h&=\frac{\mu}{2}\ms{D}_{\bs \theta^i}
    \ms{D}_{\bs \theta^j}I^h_1-\mu\left(-\frac{1}{(J^h)^2}\ms D_{\bs
    \theta^i}J^h \otimes \ms D_{\bs \theta^j}J^h+\frac{1}{J}\ms{D}_{\bs
    \theta^i} \ms{D}_{\bs \theta^j}J^h\right)\\
    &+\frac{\lambda}{4}\left[2\ms D_{\bs \theta^i}J^h \otimes
    \ms{D}_{\theta^i}J^h+2J^h \ms{D}_{\theta^i}\ms{D}_{\theta^j}J^h-2\left(
    -\frac{1}{(J^h)^2} \ms D_{\bs \theta^i}J^h\otimes \ms
    D_{\bs\theta^j}J^h+\frac{1}{J}\ms{D}_{\bs\theta^i}\ms D_{\bs\theta^j}J^h\right)\right]
\end{align}
The deformed shapes of the ring for loads of 4, 8  and 12 $kN$ are shown in
Fig.\ref{fig:srDefConfig} for different finite element meshes. From the deformed
shapes, we see that the ring has undergone extremely large deformation,
an evidence for the efficacy of our modelling-cum-discretization approach. The convergence plots of
deformation 1-forms and first Piola stress are in Fig.
\ref{fig:srConvField}. 

\section{Conclusions}
In this article, we have discussed an advanced discretization scheme for non-linear
elasticity by bringing together ideas from geometry and finite dimensional
approximations of differential forms. The numerical technique is based on a mixed
variational principle, whose input arguments are differential forms. We have applied
our discretization scheme to certain benchmark problems and, in the process, assessed its
performance against shear and volume locking as well as hour-glassing. The numerical study confirms
that our FE formulation is free of such numerical difficulties. In the literature,
almost all numerical methods which avoid these numerical instabilities use some
form of stabilization.  The absence of any such stablization term in our
variational formulation highlights the superiority of the proposed formulation as well as the
discretization technique.

An interesting extension of the present method would be to include
inertia in the formulation, wherein schemes based on variational integrators may be used to develop approximations that 
conserve energy and momenta.  An important aspect of
Cartan's moving frame that was not explored in our variational principle is the
evolution of the affine connection. For a three dimensional elastic body, such
an evolution is redundant since the Euclidean hypothesis on the geometry of the deformed
configuration completely fixes it. However, for an elastic shell, the scenario is
different. Here the affine connection has to be sought as a solution along with the
deformed configuration, when formulated as a mixed problem. We believe that constructing variational 
principles that take the
affine connection as an input is an important direction of future research.
These variational principles may then be exploited to build numerical
techniques for problems in shell theory and defect mechanics.  An equally
important research direction is the inclusion of other kinds of
interactions like electromagnetism. For the case of electromagnetism, the
present formulation is advantageous since electromagnetic fields may also be
described as differential forms. 
\section*{Acknowledgements}
BD and JK were supported by ISRO through the Centre of Excellence in Advanced
Mechanics of Materials grant No. ISRO/DR/0133. 
\appendix
\section{Discrete approximation of pulled-back area-forms}
The discrete HW principle presented in \eqref{eq:discreteHW} requires  the 
 pulled-back area-forms for their description.
Here we present a discretization for the pulled-back area-form in terms of the discrete deformation $1-$ forms given in \eqref{eq:approxTheta}. Formally, the discrete approximation
to the pulled-back area forms may be written as,
\begin{equation}
     \mathsf A^1_h=\theta^2_h\wedge \theta^3_h;\quad
     \mathsf A^2_h=\theta^3_h\wedge \theta^1_h;\quad
     \mathsf A^3_h=\theta^1_h\wedge \theta^2_h.
    \label{eq:fdAreaForms}
\end{equation}
The anti-symmetry of the area-forms along with the finite dimensionality of the approximation permits us to find a matrix representation
for the discrete area-forms. However, the entries of this matrix are not real
numbers but 2-form. Using a finite dimensional approximation
\eqref{eq:approxTheta} for $\theta^i$  in \eqref{eq:fdAreaForms}, we have,
\begin{equation}
 \theta^i\wedge\theta^j=
\begin{bmatrix}
    \bar{\theta}^i_{1} \\
    \bar{\theta}^i_{2} \\ 
    \vdots\\
    \bar{\theta}^i_{m}
\end{bmatrix}^t 
    \begin{bmatrix} 
        0 & \phi^{1}\wedge\phi^{2} & \hdots & \phi^{1}\wedge \phi^{m} \\
        -\phi^{1}\wedge\phi^{2} &0 &\hdots & \phi^{2}\wedge \phi^{m} \\ 
        \vdots & \vdots & \ddots & \vdots \\
        -\phi^{1}\wedge\phi^{m} & -\phi^{2}\wedge \phi^{m} & \hdots& 0 
    \end{bmatrix}
    \begin{bmatrix}
        \bar{\theta}^j_{1} \\
        \bar{\theta}^j_{2} \\ 
        \vdots\\
        \bar{\theta}^j_{m}
    \end{bmatrix}
    \label{eq:discretWedgeProduct}
\end{equation}
The above relationship can be formally written as $ \ms A^k_h=(\bs\theta^i)^t
\tb{W} \bs\theta^j$, where the matrix $\tb{W}$ is skew-symmetric. The matrix $\tb{W}$  will  retain its skew-symmetry as long as the
FE spaces used to approximate $\theta^i$ and $\theta^j$ are the same. The definition
of discrete area-form given in \eqref{eq:discretWedgeProduct} can be understood
as a discrete wedge operator between two 1-forms. The wedge operator is now given
by a skew symmetric matrix $\tb W$, which depends on the approximation space for
$\theta^i$. Using the definition of the discrete wedge operator, it is now easy
to compute the first derivative of discrete area form which can be written as,
\begin{equation}
    \ms D_{\bs\theta^i}  \mathsf A^k=\tb W \bs \theta^j; \quad i \neq j \neq k
\end{equation}
Similarly, the second derivative of the area forms can be written as,
\begin{equation}
    \ms D_{\bs\theta^j \bs\theta^i} \mathsf A^k=\tb W; \quad i \neq j \neq k
\end{equation}

\section{Discrete approximation of stored energy density function}
In Section 4, we introduced the discrete approximation to the stored energy
function. In the numerical simulations, we assumed that the stored energy
function depends on the invariants of the pulled-back metric
tensor $C$. We first discretize the invariants of $C$. A discrete approximation to the stored energy function may be constructed by taking the discrete invariants as the input arguments. The discrete
stored energy function can thus be written as,
\begin{equation}
    W^h:=W(I_1^h, I^h_2, J^h).
\end{equation}
In the above definition, $I_1^h$, $I_2^h$ and $J^h$ denote the discrete
approximations for the invariants of the pulled-back metric tensor.

\subsection{Discrete approximation of $I_1$ and its derivatives}
Using the definition of $I_1$ given in \eqref{eq:I1}, its discrete approximation
may be written as,
\begin{equation}
    I_1^h:=\sum_{i=1}^3\left\langle\theta^i_h,\theta^i_h\right\rangle
\end{equation}
The right hand side of the above definition is a real valued bilinear form. In terms of the shape functions, this bilinear form can be
written as,
\begin{equation}
    I_1^h:=\sum_{i=1}^3 (\bs \theta^i)^t \bs \phi ^t \bs \phi \bs \theta^i
\end{equation}
The first derivative of $I_1^h$ with respect to the DoF vector may now be
computed as,
\begin{equation}
    \ms D_{\bs \theta^i} I_1^h=2 \bs \phi^t \bs \phi \bs \theta^i 
\end{equation}
The second derivative with respect to the DoF's may be evaluated as,
\begin{equation}
    \ms D_{\bs \theta^j \bs \theta^i}I_1=
    \left\{
        \begin{array}{ll}
             2\bs \phi^t \bs \phi;&\quad i=j\\
             \tb 0;&\quad i\neq j
        \end{array}
    \right.
\end{equation}
\subsection{Discrete approximation of $I_{2}$ and its derivatives}
The discrete approximation of the second invariant of $C$ may be written as,
\begin{equation}
    I_2^h:=\sum_{i=1}^3\left\langle  \ms A^i_h, \ms A^i_h \right\rangle
\end{equation}
where the discrete approximations for $ \ms A^i$ are given in
\eqref{eq:fdAreaForms}. 
\begin{align}
    \ms D_{\bs \theta^1}I_2^h&=-2(\bs \phi \wedge \theta^3_h) \ms A^2_h+2(\bs \phi
    \wedge \theta^2_h )\ms  A^3_h\\
    \ms D_{\bs \theta^2}I_2^h&=2(\bs \phi \wedge \theta^3_h)\ms A^1_h-2(\bs \phi
    \wedge \theta^1_h )\ms A^3_h\\
    \ms D_{\bs \theta^3}I_2^h&=-2(\bs \phi \wedge \theta^2_h)\ed \ms A^1_h+2(\bs \phi
    \wedge \theta^1_h )\ms  A^2_h
\end{align}
The second derivatives of $I_2$ with respect to the DoFs may be computed as,
\begin{align}
    \ms D_{\bs \theta^1 \bs \theta^1}I_2^h&=2(\bs \phi \wedge \theta^2)^t(\bs
    \phi\wedge \theta^2)+2(\bs \phi \wedge\theta^3)^t(\bs \phi \wedge
    \theta^3)\\
    \ms D_{\bs \theta^2 \bs \theta^2}I_2^h&=2(\bs \phi \wedge \theta^3)^t(\bs
    \phi\wedge \theta^3)+2(\bs \phi \wedge\theta^1)^t(\bs \phi \wedge
    \theta^1)\\
    \ms D_{\bs \theta^3 \bs \theta^3}I_2^h&=2(\bs \phi \wedge \theta^1)^t(\bs
    \phi\wedge \theta^1)+2(\bs \phi \wedge\theta^2)^t(\bs \phi \wedge
    \theta^2)\\
    \ms D_{\bs \theta^1 \bs \theta^2}I_2^h&=-2[\bs \phi^t(\bs
    \phi\wedge \ms A^3)+(\bs \phi \wedge\theta^1)^t(\bs \phi \wedge
    \theta^2)]\\
    \ms D_{\bs \theta^2 \bs \theta^3}I_2^h&=-2[\bs \phi^t(\bs
    \phi\wedge \ms A^1)+(\bs \phi \wedge\theta^2)^t(\bs \phi \wedge
    \theta^3)]\\
    \ms D_{\bs \theta^3 \bs \theta^2}I_2^h&=-2[\bs \phi^t(\bs
    \phi\wedge \ms A^2)+(\bs \phi \wedge\theta^3)^t(\bs \phi \wedge
    \theta^1)]
\end{align}

\subsection{Discrete approximation of $J$ and its derivatives}
Using the definition of $J$ and FE approximation of deformation 1-forms, the
discrete approximation for $J$ can be written as,
\begin{equation}
    \nonumber
    J^h:=\theta^1_h\wedge\theta^2_h\wedge\theta^3_h.
    \label{ed:disJ}
\end{equation}
The relationship between the deformation 1-form DoFs and $J^h$ is trilinear and anti-symmetric. 
Using the definition of $\theta^1_h$ in the above expression, we
have, 
\begin{equation}
    J^h=(\bs \phi \bs \theta^1)\wedge(\bs \phi \bs \theta^2)\wedge(\bs \phi \bs
    \theta^3).
\end{equation}
The first derivative can thus be computed as,
\begin{equation}
    \ms D_{\theta^k}J^h=-1^{k+1}\bs\phi\wedge \theta^i_h\wedge\theta^j_h;\quad i\neq
    j\neq k; \quad i<j.
\end{equation}
Similarly, the second derivative with respect to the DoFs can be computed as,
\begin{align}
    \nonumber
    \ms D_{\bs \theta^1 \bs \theta^2}J^h&=\bs \phi \wedge \bs \phi \wedge
    \theta^3\\
    \ms D_{\bs \theta^1 \bs \theta^3}J^h&=-\bs \phi \wedge \bs \phi \wedge
    \theta^2\\
    \nonumber
    \ms D_{\bs \theta^2 \bs \theta^3}J^h&=-\bs \phi \wedge \bs \phi \wedge
    \theta^1,
\end{align}
We also have $\ms D_{\bs \theta^i \bs \theta^i}J^h=\tb 0$.

\bibliographystyle{abbrv}
\bibliography{mfem}
\end{document}